\documentclass[final,3p,times]{elsarticle}
\usepackage{url,microtype,subcaption}
\usepackage[onehalfspacing]{setspace}
\usepackage[thinc]{esdiff} 
\usepackage{amssymb}
\usepackage{amsmath}
\usepackage{bm}
\usepackage{tikz}
\usepackage{upgreek}
\usepackage{adjustbox}
\usepackage{float} 
\usepackage{adjustbox}
\usetikzlibrary{patterns}
\usetikzlibrary{arrows}
\usepackage{url}
\usepackage[switch,columnwise]{lineno}
\usepackage[colorlinks=true,linkcolor=blue,anchorcolor=blue,citecolor=blue,filecolor=blue,menucolor=blue,runcolor=blue,urlcolor=blue,breaklinks]{hyperref}
\definecolor{PBMcolor}{rgb}{1,0,0}
\definecolor{DDMcolor}{rgb}{0.9,0.6,0}
\definecolor{CoSTAcolor}{rgb}{0,0,1}
\definecolor{PGNNCoSTAcolor}{rgb}{0,0.5,0}\definecolor{EintXc}{rgb}{.7,0,0}
\definecolor{EintIc}{rgb}{0,0.6,0.2}

\DeclareRobustCommand{\pbmline}{\raisebox{1mm}{\tikz{
    \draw[-,PBMcolor,solid,line width = 1.0pt](0.,0mm) -- (5mm,0mm);
}}}

\DeclareRobustCommand{\ddmerr}{\tikz{
    \fill[DDMcolor!40!white,solid,line width = 1.0pt](0.mm,0) rectangle (5mm,2mm);
    \draw[-,DDMcolor,solid,line width = 1.0pt](0.,0mm) -- (5mm,0mm);
}}
\DeclareRobustCommand{\CoSTAerr}{\tikz{
    \fill[CoSTAcolor!40!white,solid,line width = 1.0pt](0.mm,0) rectangle (5mm,2mm);
    \draw[-,CoSTAcolor,solid,line width = 1.0pt](0.,0mm) -- (5mm,0mm);
}}

\DeclareRobustCommand{\disc}[1]{
\begin{tikzpicture}
\fill[#1] (0,0) circle (1.2 mm);
\end{tikzpicture}
}

\newcommand{\subfigurecaptionInterpol}[1]{\caption{Solution \textit{#1} with $\alpha\!=\!0.7$ (left) and $\alpha\!=\!1.5$ (right)}}
\newcommand{\subfigurecaptionExtrapol}[1]{\caption{Solution \textit{#1} with $\alpha\!=\!-0.5$ (left) and $\alpha\!=\!2.5$ (right)}}
\usepackage{pgfplots}
\pgfplotsset{compat=1.17}

\makeatletter
\def\addlegendimage{\pgfplots@addlegendimage}

\usepackage{pifont}
\usepackage[ruled,vlined]{algorithm2e}
\usepackage{algorithmic}
\usepackage{booktabs} 
\usepackage{siunitx}
\modulolinenumbers[5]

\newbox\tempbox%
\newenvironment{nomenclature}{\par\vskip6pt plus1pt\setbox\tempbox\vbox\bgroup\if@twocolumn\hsize17.5pc\else\hsize37.3pc\fi\vspace*{-17.7pt}\section*{Nomenclature}}{\egroup\par\noindent\fboxsep10pt\fbox{\box\tempbox}\vskip6pt}
\def\deflist{\@ifnextchar[{\@deflist}{\@deflist[]}}
\newbox\defbox%
\newdimen\defboxdimen%
\def\@deflist[#1]{%
\setbox\defbox\hbox{#1:\quad}%
\defboxdimen\wd\defbox%
\def\tempa{#1}%
\par\addvspace{12pt plus2pt minus2pt}%
\setlength\parindent{0pt}%
\def\deftitle##1{{\noindent\itshape##1}\par}
\ifx\tempa\@empty%
\def\defitem##1{\@hangfrom{##1\ }}
\else%
\def\defitem##1{\@hangfrom{\hbox to \defboxdimen{##1\quad}}}
\fi%
\def\defterm##1{##1\par}
}

\usepackage{soul}

\journal{Journal of \LaTeX\ Templates}

\begin{document}

\begin{frontmatter}
    \title{Enhancing Elasticity Models: A Novel Corrective Source Term Approach for Accurate Predictions}
    
    \author[trondsaddress]{Sondre Sørbø} \ead{sondre.sorbo@sintef.no}
    \author[trondsaddress,sindresaddress]{Sindre Stenen Blakseth} \ead{sindre.blakseth@sintef.no}
    \author[adilsaddress,SINTEFaddress]{Adil Rasheed\corref{mycorrespondingauthor}}
    \cortext[mycorrespondingauthor]{Adil Rasheed}
    \ead{adil.rasheed@ntnu.no}
    \author[trondsaddress,SINTEFaddress]{Trond Kvamsdal} \ead{trond.kvamsdal@ntnu.no}
    \author[omersaddress]{Omer San} \ead{osan@utk.edu }

    \address[trondsaddress]{Department of Mathematical Sciences, Norwegian University of Science and Technology}
    \address[adilsaddress]{Department of Engineering Cybernetics, Norwegian University of Science and Technology}
    \address[SINTEFaddress]{Mathematics and Cybernetics, SINTEF Digital}
    \address[omersaddress]{Department of Mechanical, Aerospace and Biomedical Engineering, University of Tennessee, Knoxville}
    \address[sindresaddress]{Department of Gas Technology, SINTEF Energy Research}

\begin{abstract}
\textit{Motivation:} In the context of safety-critical applications amid the recent wave of digitalization, there is a pressing demand for computational models that are simultaneously efficient, accurate, generalizable, and trustworthy. While physics-based models have long been employed to simulate complex phenomena and are known for their trustworthiness and generalizability, they often fall short when real-time performance is essential. The simplifications required to make them computationally viable, such as reducing problem dimensions or linearizing non-linear characteristics, can compromise their accuracy. On the other hand, data-driven models offer computational efficiency and accuracy but often lack generalizability and operate opaquely, making them less suitable for safety-critical applications. This research aims to bridge these gaps by employing a data-driven approach to correct simplifications in physics-based models.

\textit{Method:} This article introduces a novel methodology referred to as the "corrective source term approach" to address the limitations of simplified physics-based models. This approach leverages data-driven techniques to rectify two common simplifications: dimension reduction (e.g., from 3D to 2D) and linearization of non-linear problem characteristics. By integrating data-driven corrections into the physics-based framework, the research demonstrates a hybrid modeling approach that enhances the accuracy and reliability of physics-based models.

\textit{Main Result:} To assess the effectiveness of the proposed methodology, it is applied to model various elasticity problems. The results of this study showcase the superior performance of the corrective source term approach compared to traditional physics-based models and end-to-end data-driven models. Specifically, the hybrid model exhibits significant improvements in terms of accuracy, reduced model uncertainty, and enhanced generalizability, making it a promising solution for safety-critical applications where computational efficiency, accuracy, and trustworthiness are paramount concerns.
\end{abstract}

\end{frontmatter}

\begin{nomenclature}
\begin{deflist}[AAAAA] 
\defitem{PBM}\defterm{Physics-Based Modeling}
\defitem{ML}\defterm{Machine Learning}
\defitem{DDM}\defterm{Data-Driven Modeling}
\defitem{HAM}\defterm{Hybrid Analysis and Modeling}
\defitem{NN}\defterm{Neural Network}
\defitem{RNN}\defterm{Recurrent Neural Network}
\defitem{ROM}\defterm{Reduced-Order Modeling}
\defitem{ODE}\defterm{Ordinary Differential Equation}
\defitem{DNN}\defterm{Deep Neural Network}
\defitem{PINN}\defterm{Physics-Informed Neural Network}
\defitem{PDE}\defterm{Partial Differential Equation}
\defitem{PGML}\defterm{Physics-Guided Machine Learning}
\defitem{CoSTA}\defterm{Corrective Source Term Approach}
\defitem{$\boldsymbol{u}$}\defterm{Displacement}
\defitem{$\boldsymbol{\sigma}$}\defterm{Cauchy stress tensor}
\defitem{$n$}\defterm{Number of spatial dimensions}
\defitem{$\boldsymbol{f}$}\defterm{Imposed structural load}
\defitem{$t$}\defterm{Time}
\defitem{$\boldsymbol{\epsilon}$}\defterm{Strain}
\defitem{$\boldsymbol{C}$}\defterm{Constitutive tensor relating stress and strain}
\defitem{$E$}\defterm{Young's modulus}
\defitem{$\nu$}\defterm{Poisson ratio}
\defitem{$\mathbf{D}$}\defterm{Strain operator}
\defitem{$\gamma$}\defterm{Engineering shear strain}
\defitem{$\Omega$}\defterm{Spatial domain}
\defitem{$T$}\defterm{Final time}
\defitem{IBVP}\defterm{Initial Boundary Value Problem}
\defitem{$\boldsymbol{x}$}\defterm{Position}
\defitem{$\partial\Omega$}\defterm{Domain boundary}
\defitem{$\partial\Omega_{\textsc{d}}$}\defterm{Domain boundary segment with Dirichlet boundary conditions}
\defitem{$\partial\Omega_{\textsc{n}}$}\defterm{Domain boundary segment with Neumann boundary conditions}
\defitem{$\boldsymbol{n}$}\defterm{Boundary normal vector}
\defitem{$g_{\textsc{d}}$}\defterm{Dirichlet boundary data}
\defitem{$g_{\textsc{n}}$}\defterm{Neumann boundary data}
\end{deflist}
\end{nomenclature}
\begin{nomenclature}
\begin{deflist}[AAAAA] 
\defitem{$\boldsymbol{\mu}_0$}\defterm{Initial condition for displacement}
\defitem{$\boldsymbol{\rho}_0$}\defterm{Initial condition for the time derivative of displacement}
\defitem{VP}\defterm{Variational Problem}
\defitem{X}\defterm{$n$-dimentional Hilbert space on the spacial domain ($\Omega$)}
\defitem{$X_D$}\defterm{Subset of $X$ satisfying boundary conditions}
\defitem{$X_0$}\defterm{Subspace of $X$ vanishing at the boundary}
\defitem{$Y(\cdot)$}\defterm{Solution search space with spacial search space $\cdot$}
\defitem{$\boldsymbol{u}_{\textsc{d}}$}\defterm{Lifting function for non-homogeneous boundary conditions}
\defitem{$\boldsymbol{u}_{\textsc{0}}$}\defterm{Unknown part of displacement with homogeneous boundary conditions}
\defitem{$\boldsymbol{v}$}\defterm{Test function}
\defitem{$\boldsymbol{w}$}\defterm{Test function}
\defitem{$a(\cdot, \cdot$)}\defterm{Bilinear form from the variational problem formulation}
\defitem{$(\cdot, \cdot$)}\defterm{Bilinear form from the variational problem formulation}
\defitem{$l(\cdot$)}\defterm{Linear form from the variational problem formulation}

\defitem{$\cdot_h$}\defterm{Discrete approximation of $\cdot$}
\defitem{$N_{\mathrm{dof}}$}\defterm{Number of degree of freedom in finite element approximation}
\defitem{$N_{\mathrm{el}}$}\defterm{Number of elements in the domain triangulation}
\defitem{$\phi_i$}\defterm{Nodal basis function}
\defitem{$({u}_{h,0})_j(t)$}\defterm{Coefficient $j$ of the nodal basis representation of the trial function}
\defitem{${v}_i$}\defterm{Coefficient $i$ of the nodal basis representation of the test function}
\defitem{$P_1(\cdot)$}\defterm{The space of linear polynomials on $\cdot$}
\defitem{$A$}\defterm{Stiffness matrix}
\defitem{$F$}\defterm{Load vector}
\defitem{$M$}\defterm{Mass matrix}
\defitem{$U$}\defterm{Finite element approximation of displacement}
\defitem{$\cdot^{(i)}$}\defterm{Quantity $\cdot$ evaluated at time level $i$}
\defitem{$k$}\defterm{Time step}
\defitem{$N$}\defterm{Number of neural network layers}
\defitem{$D$}\defterm{Mapping defined by a neural network}
\defitem{$\mathcal{T}_i$}\defterm{Affine transformation in layer $i$ of a neural network}
\defitem{$s_i$}\defterm{Activation function in layer $i$ of a neural network}
\defitem{$d_i$}\defterm{Number of nodes in layer $i$ of a neural network}
\defitem{$L_{\Omega}$}\defterm{General differential operator}
\defitem{$L_{\partial\Omega}$}\defterm{General differential operator}
\defitem{$\omega$}\defterm{Unknown of interest in general differential equation}
\defitem{$\theta$}\defterm{General source term}
\defitem{$\psi$}\defterm{General boundary condition}

\end{deflist}
\end{nomenclature}
\begin{nomenclature}
\begin{deflist}[AAAAA] 
\defitem{$\tilde{\cdot}$}\defterm{Perturbed instance of $\cdot$ (perturbation caused by e.g.\ some error or noise)}
\defitem{$r$}\defterm{Residual between perturbed and unperturbed differential equation}
\defitem{$\tilde{\tilde{\cdot}}$}\defterm{Corrected instance of perturbed quantity $\tilde{\cdot}$}
\defitem{$\alpha$}\defterm{Parametrization of system state}

\defitem{$\mathcal{A}_{\mathrm{train}}$}\defterm{Set of $\alpha$-values used for training}
\defitem{$\mathcal{A}_{\mathrm{val}}$}\defterm{Set of $\alpha$-values used for validation}
\defitem{$\mathcal{A}_{\mathrm{test}}$}\defterm{Set of $\alpha$-values used for testing}
\defitem{$K$}\defterm{Number of time steps}
\defitem{$\check{\cdot}$}\defterm{Vector form of exact solution}
\defitem{$\bar{\cdot}$}\defterm{Prediction based on previous exact step}
\defitem{$\hat{\cdot}$}\defterm{Prediction based on previous predicted step}
\defitem{$x$}\defterm{First spacial coordinate}
\defitem{$y$}\defterm{Second spacial coordinate}
\defitem{$z$}\defterm{Third spacial coordinate}
\defitem{RRMSE}\defterm{Relative Root Mean Square Error}
\end{deflist}
\end{nomenclature}

\section{Introduction} \label{sec:introduction}
Predictive modeling and simulation has traditionally been dominated by physics-based modeling (PBM). With the rise of machine learning (ML) in recent times, data-driven modeling (DDM)\footnote{Throughout the text, we will use the acronym DDM to refer to both data-driven modeling and data-driven models. Similarly, PBM may refer to physics-based modeling or physics-based models.} has shown its ability to outperform PBM in many situations~\citep{Wei2018PredictingTE, Baldi2016JetSC, IbarraBerastegi2015ShorttermFO, Jia2021PhysicsGuidedML}. However, DDM comes with its own disadvantages, limiting these methods' overall usefulness. In \cite{san2021hybrid}, the authors describe the ideal model in the context of digital twins as generalizable, trustworthy, computationally efficient, accurate and self-evolving. A model's generalizability is its ability to solve various problems without problem-specific fine-tuning. Trustworthiness refers to the extent to which a model is explainable or interpretable, while computational efficiency and accuracy refer to the model's ability to make real-time predictions that match ground truth as closely as possible.
Lastly, a model is self-evolving if it can learn and evolve when new situations are encountered. PBMs, when based on the correct physics, can be accurate and generalizable but are usually computationally demanding and do not adapt to new scenarios automatically. DDMs, on the other hand, after training, are very efficient, possibly very accurate and can be self-evolving. However, they typically lack in trustworthiness and generalizability.

For both PBM and DDM, a model's accuracy depends on the knowledge used to build the model. To achieve high predictive accuracy, PBMs require a proper understanding of all relevant physical phenomena, as well as mathematical methods for solving the equations used to represent these phenomena. If our physics knowledge is inaccurate and/or we lack methods for efficient computation, PBMs' accuracy will be significantly reduced.
Contrary to PBMs, DDMs do not rely on physics knowledge in order to achieve high accuracy. Instead, they rely on (possibly large amounts of) data that is representative for the scenarios in which the model is to be used. With good data, DDMs can deliver highly accurate results at only a fraction of the computational cost of PBMs with comparable accuracy, if such a PBM is even available. However, given unrepresentative data, DDMs perform poorly due to their lack of generalizability. Moreover, since DDMs rely on data rather than explainable and verifiable physics knowledge, their interpretability is lacking. This combination of poor generalizability and interpretability results in poor trustworthiness and has prevented widespread utilization of DDM in high-stakes applications where even a single poor prediction can be detrimental.

The steadily increasing industrial adoption of digital twins \citep{8972429} entails that models possessing all the characteristics identified by \cite{san2021hybrid} are as relevant as ever. However, as we have seen, neither PBM nor DDM possess all four of these characteristics. As such, their respective deficiencies imply that neither paradigm is suitable for reaping the full benefits of the real-time monitoring and control enabled by digital twins.

To counter these deficiencies, hybrid analysis and modeling (HAM) is emerging as a new eclectic paradigm that combines techniques from both PBM and DDM. The HAM approach unites the advantages of PBM, such as generalizability, interpretability, solid foundation, and comprehension, with the accuracy, computational efficiency, and automatic pattern-recognition abilities of DDM. In their recent reviews, \cite{willard2020integrating} and \cite{san2021hybrid} offer a comprehensive look at techniques for combining DDM with PBM. A lot of the hybridization techniques can be classified into the following categories:

(i) Embedding PBMs inside NNs, 
(ii) Model order reduction, 
(iii) Physics-based regularisation terms, 
(iv) Data-driven equation discovery, 
(v) Error correction approaches, 
(vi) Sanity check mechanisms using PBMs.

In the following sections, we present related work and discuss the advantages and disadvantages of the approaches.

    \subsection{Methods for embedding PBMs directly into NNs}
    This approach to hybridisation is the most straightforward. More advanced techniques usually create a differentiable PBM that can be used as a layer in a neural network. For instance, OptNet \cite{amos_optnet_2017} is a differentiable convex optimisation solver that can be used as a layer in a network. In \cite{avilabelbuteperes_end_2018} the authors proposed the differentiable physics engine, a rigid body simulator that can be embedded in a NN. They showed that it is possible to learn a mapping from visual data to the positions and velocities of objects, which are then updated using the simulator. Authors in \cite{yu2020sds} simulated a structural dynamics problem by designing a hybrid recurrent NN (RNN) that contains an implicit numerical integrator. These approaches are usually quite data-efficient, but they can make both inference and training more expensive.

    \subsection{Model order reduction methods}
    The reduced-order modelling (ROM) approach has been widely used to project complex partial differential equations onto a reduced dimensional space based on the singular value decomposition of the offline high fidelity simulation data, resulting in a set of ordinary differential equations (ODEs) which are much faster to solve \citep{quarteroni2014reduced,fonn2019dcp}. This has enabled high-fidelity numerical solvers to be accelerated by several orders of magnitude. However, ROMs can become unstable in the presence of unknown/unresolved complex physics. To address this issue, recent research has demonstrated how unknown and hidden physics within a ROM framework can be accounted for using deep neural networks (DNNs) \citep{pawar2019aet,pawar2020ddr}. Nevertheless, ROMs require the exact form of the original equation before they can be applied. 

    \subsection{Physics-based regularization terms} By incorporating a physics-based model (PBM) into the objective function, deep learning models (DDMs) can be guided to adhere to known physical laws during training. An example of this is the physics-informed neural network (PINN) proposed by \cite{raissi2019physics}, which uses a neural network (NN) to represent the solution to a partial differential equation (PDE) and adds an additional loss term to penalize any deviations from the equation at a sample of points. \cite{zobeiry_physics_2021} applied the PINN approach to solve heat transfer problems in manufacturing processes, while \cite{arnold_state_2021} extended it to enable control in a state-space setting. \cite{shen_physics_2021} created a model to classify the health of the bearings by training a NN on physics-based features and regularizing the model using the output from a physics-based threshold model. More applications of PINN can be found in heat transfer modeling \cite{BILLAH2023106336}, multicomponent reactor modeling \citep{SUN2023103525} or in predicting the lifetime under multiaxial loading \cite{HALAMKA2023109351}. However, these approaches require precise knowledge of the loss term and can be difficult to train due to the increased complexity of the cost function, particularly if the regularization term necessitates the evaluation of a complex model. 
    
    \subsection{Data-driven equation discovery} 
    Sparse regression, which is based on $l_1$ regularization, and symbolic regression, which is based on gene expression programming, have been demonstrated to be highly successful in uncovering hidden or partially known physical laws from data. Examples of this type of approach can be found in \cite{Champion22445} and \cite{vaddireddy2020fes}. On a more applied side, the authors in \cite{ZHANG2023112349} discovered the equation describing the relationship between features and material characteristics while Meyer et al \cite{MEYER2023105416} utilized thermodynamically consistent neural network to model plasticity and discover the associated evolution laws. However, there are some drawbacks to this class of methods. For instance, with sparse regression, extra features must be manually created, while with symbolic regression, the resulting models can be unstable and prone to overfitting.

    \subsection{Physics guided machine learning}
     On many occasions part of the physics governing a process is known, but the actual form describing is not known. To exploit such partial knowledge, Pawar et al \cite{pawar2021pgml} proposed a physics-guided machine learning (PGML) approach. The basic idea behind the PGML approach is to inject partial knowledge into one of the layers within a DNN to guide the training process. Partially knowledge can, for example, come from a simplistic model as has been shown in \cite{pawar2021pgml}, \cite{Robinson2022pgn}. More recently, the Theseus library \cite{pineda2022theseus} provides a framework for conducting guided training of neural networks. However, this approach does not take advantage of the partially known form of the equation.

    \subsection{Proposal}
    From the previous discussion, it is clear that almost all HAM approaches have pros and cons. One limitation of all the methods is that they are more tilted towards the data-driven modeling approach. To this end, Corrective Source Term Approach (CoSTA) has been proposed recently. CoSTA is a method proposed by \cite{Blakseth2022dnn} that explicitly addresses the problem of unknown physics. This is done by augmenting the governing equations of a PBM describing partial physics with a DNN-generated corrective source term that takes into account the remaining unknown/ignored physics. One added benefit of the CoSTA approach is that the physical laws can be used to keep a sanity check on the predictions of the DNN used, i.e. checking conservation laws. A similar approach has also been used to model unresolved physics in turbulent flows \citep{maulik2019sgm,pawar2020apa}. The method has also been shown to work well for modeling complex aluminium extraction process \cite{robinson2022novel}. In this regard, the main contributions distinguishing the present work from previous publications on CoSTA \citep{blakseth2021ica,Blakseth2022cpb,Blakseth2022dnn, robinson2022novel}, are summarized below.
\begin{itemize}
    \item We investigate whether CoSTA can be used to correct modeling errors incurred due to \textbf{dimensionality reduction} in PBMs. Reducing the dimensionality of a model is commonly seen in engineering applications in order to reduce computational complexity, thereby achieving real-time performance. However, this comes at the cost of reduced predictive accuracy. CoSTA has not previously been used to correct for this kind of modeling error.
    \item We investigate whether CoSTA can be used to correct modeling errors incurred due to \textbf{linearization} of non-linear governing equations. Linearization is another technique that is commonly used to speed up models, and which also reduces predictive accuracy.
    The use of CoSTA for correcting linearization error has been touched upon briefly in \citep{blakseth2021ica}, but is considered in much greater detail in the present work.
    \item We study the effect of randomness in DNN initialization and training procedures on CoSTA's performance, relative to stand-alone PBM and DDM.
    \item We demonstrate how to combine CoSTA with PBMs based on the \textbf{finite element method}. Previously, CoSTA has only been used in conjunction with finite volume methods.
    \item We apply CoSTA to a new class of problems: \textbf{elasticity modeling}. Previously, CoSTA has only been used to model heat diffusion, which is a fundamentally different phenomenon.
\end{itemize}

The article is structured as follows: Section~\ref{sec:theory} presents the relevant theory for the work. Section~\ref{sec:methodology} presents the method applied in the paper and a description of the cases considered. The results are presented and discussed in section~\ref{sec:results}. Finally, in Section~\ref{sec:conclusionandfuturework}, conclusions are given, and potential future work is presented.     

\section{Theory}\label{sec:theory}
In this section, we discuss the theory underlying the models used in our numerical experiments. We begin with general considerations regarding PBM in Section~\ref{subsec:pbm}, while Sections~\ref{subsec:lin_elasticity} and~\ref{sec:nonlin_elasticity} are devoted to the physics-based elasticity models used in the present work. The latter sections are largely based on the textbooks by \cite{linear_elasticity} and \cite{irgens2008continuum}. DDM (Section~\ref{subsec:ddm}), CoSTA for elasticity problems (Section~\ref{subsec:costa}) and the method of manufactured solutions (Section~\ref{subsec:manufacturedsolution}) are the remaining topics covered in this section.

\subsection{Physics Based Modeling}\label{subsec:pbm}

PBM involves careful observation of a physical phenomenon of interest (elasticity in the current work), development of its partial understanding, expression of the understanding in the form of mathematical equations and ultimately solution of these equations. Due to the partial understanding and numerous assumptions along the steps from observation to solution of the equations, a large portion of the important governing physics gets ignored. Due to high computational costs, high fidelity simulators with minimal assumptions have so far been limited to the offline design phase only. Despite this major drawback, what makes these models attractive are sound foundations from first principles, interpretability, generalizability and existence of robust theories for the analysis of stability and uncertainty. The PBMs used in this paper are based on partial differential equations (PDEs) describing linear elasticity in solid materials (cf.\ Section~\ref{subsec:lin_elasticity}). These PDEs are discretized using the finite element method\footnote{See e.g.\ \citep{quarteroni2017nummod} or \citep{brenner2008FEM} for excellent introductions to finite element method.} along the spatial dimension, and the backward Euler method along the temporal dimension.
The numerical methods will, in the vast majority of interesting scenarios, introduce some discretization error into these models. Moreover, modeling error may be present as well due to the governing equations being simplified, incomplete or even incorrect. Sometimes, modeling error may be purposefully introduced to simplify the model, thereby increasing its computational efficiency. Other times, modeling error could stem from a lack of knowledge about the system to be modelled. In this work, we explore the impact of modeling error through PDE linearization, dimensionality reduction and unknown load terms, as described in Section~\ref{sec:methodology}.

\subsection{Linear Elasticity Modeling}
\label{subsec:lin_elasticity}
The purpose of elastic modeling is to predict the response of a solid material to external, typically mechanical loads. This response is quantified in terms of the displacement $\boldsymbol{u}$ and stresses~$\boldsymbol{\sigma}$. For an $n$-dimensional system, $\boldsymbol{u}$ is a $n$-dimensional vector field whose elements describe a point's actual position in relation to where that point would have been if no forces were applied to the system. Moreover, $\boldsymbol{\sigma}$, known as the Cauchy stress tensor, is a second order tensor with dimension $n\times n$, describing the system's internal forces. The displacement and stresses are related through Newton's 2nd law written for elastic continua. For isotropic systems with unit mass density, this form of Newton's 2nd law reads
\begin{equation}
    \nabla \cdot \boldsymbol{\sigma} - \ddot{\boldsymbol{u}} = - \boldsymbol{f}.
    \label{eq:le1}
\end{equation}
Here, $\boldsymbol{f}$ denotes the forces (loads) acting on the system and $\ddot{\boldsymbol{u}} = \partial^2 \boldsymbol{u} / \partial t^2$ is the acceleration field where $t$ is time. Since we have two unknown fields ($\boldsymbol{\sigma}$ and $\boldsymbol{u}$), but only one equation, we need another constraint to have a well-defined problem. To achieve this we introduce a kinematic relationship between the strain field $\boldsymbol{\varepsilon}$ (that describes the relative stretching (deformation) of an infinitesimal element of the system) and the displacement $\boldsymbol{u}$, see Equation~\eqref{eq:le2}, and then a constitutive relation between the the stresses $\boldsymbol{\sigma}$ and $\boldsymbol{\varepsilon}$, see Equation ~\eqref{eq:le3}, where the $\boldsymbol{C}$ is the fourth order constitutive tensor.
\begin{align}
    \varepsilon &= \frac{1}{2}(\nabla \boldsymbol{u} + (\nabla \boldsymbol{u})^T)
    \label{eq:le2}
    \\
    \boldsymbol{\sigma} &= \boldsymbol{C} \varepsilon
    \label{eq:le3}
\end{align}

In computational mechanics, it is common to introduce the Voigt notation to represent the symmetric stress and strain tensors as vectors, $\boldsymbol{\upsigma}$ and $\boldsymbol{\upvarepsilon}$ respectively, and the fourth order constitutive tensor $\boldsymbol{C}$ as a two-dimensional matrix $\mathbf{C}$.   

For a system with constant Young's modulus, $E$, and constant Poisson ratio, $\nu$, the equations~\eqref{eq:le1}--\eqref{eq:le2} can be written more conveniently as
\begin{align}
    \boldsymbol{\upvarepsilon} &= \mathbf{D} \boldsymbol{u}
    \label{eq:2dle1}
    \\
    \boldsymbol{\upsigma} &= \mathbf{C} \boldsymbol{\upvarepsilon}
    \label{eq:2dle2}
    \\
     \mathbf{D}^\textsc{t} \boldsymbol{\upsigma} - \ddot{\boldsymbol{u}} &= - \boldsymbol{f},
    \label{eq:2dle3}
\end{align}
where the precise definitions of $\mathbf{C}$, $\boldsymbol{\upsigma}$, $\boldsymbol{\upvarepsilon}$, and the strain giving differantial operator $\mathbf{D}$ depend on the number of space dimensions~$n$. We will refer to Equations~\eqref{eq:2dle1}-\eqref{eq:2dle3} as the linear elasticity equations. 

In 2D, we have
\begin{equation}
    \mathbf{D}
    = 
    \left[
    \begin{array}{cc}
         \frac{\partial}{\partial x} & 0 \\
                    0                & \frac{\partial}{\partial y}\\
         \frac{\partial}{\partial y} & \frac{\partial}{\partial x}\\
    \end{array}
    \right],
    \quad
     \boldsymbol{\upvarepsilon}
    = 
    \left[
    \begin{array}{c}
         \varepsilon_{xx}\\
         \varepsilon_{yy}\\
         \gamma_{xy}\\
    \end{array}
    \right],
    \quad
    \boldsymbol{\upsigma}
    = 
    \left[
    \begin{array}{c}
         \sigma_{xx}\\
         \sigma_{yy}\\
         \sigma_{xy}
    \end{array}
    \right],
    \label{eq:epssigmadef}
\end{equation}
and
\begin{equation}
    \mathbf{C} =\frac{E}{1-\nu^2}
    \left[
    \begin{array}{ccc}
         1&\nu&0\\
         \nu&1&0\\
         0&0&\frac{1-\nu}{2}
    \end{array}
    \right],
\end{equation}
while the 3D definitions read
\begin{equation}
    \mathbf{D}
    = 
    \left[
    \begin{array}{ccc}
         \frac{\partial}{\partial x} & 0 & 0\\
         0 & \frac{\partial}{\partial y} & 0\\
         0 & 0 & \frac{\partial}{\partial z}\\
         0 & \frac{\partial}{\partial z} & \frac{\partial}{\partial y}\\
         \frac{\partial}{\partial z} & 0 &  \frac{\partial}{\partial x}\\
         \frac{\partial}{\partial y} & \frac{\partial}{\partial x} & 0\\
    \end{array}
    \right],
    \quad
    \boldsymbol{\upvarepsilon}
    = 
    \left[
    \begin{array}{c}
         \varepsilon_{xx}\\
         \varepsilon_{yy}\\
         \varepsilon_{zz}\\
         \gamma_{yz}\\
         \gamma_{zx}\\
         \gamma_{xy}\\
    \end{array}
    \right],
    \quad
    \boldsymbol{\upsigma}
    = 
    \left[
    \begin{array}{c}
         \sigma_{xx}\\
         \sigma_{yy}\\
         \sigma_{zz}\\
         \sigma_{yz}\\
         \sigma_{zx}\\
         \sigma_{xy}\\
    \end{array}
    \right].
\end{equation}
and
\begin{equation}
    \mathbf{C} =\frac{E}{(1+\nu)(1-2\nu)}
    \left[
    \begin{array}{cccccc}
        1-\nu&\nu&\nu&0&0&0\\
        \nu&1-\nu&\nu&0&0&0\\
        \nu&\nu&1-\nu&0&0&0\\
        0&0&0&\frac{1-2\nu}{2}&0&0\\
        0&0&0&0&\frac{1-2\nu}{2}&0\\
        0&0&0&0&0&\frac{1-2\nu}{2}\\
    \end{array}
    \right].
\end{equation}

Notice that, $\gamma_{yz}= 2 \varepsilon_{yz}$, $\gamma_{zx}= 2 \varepsilon_{zx}$, and $\gamma_{xy}= 2 \varepsilon_{xy}$, are the engineering shear strains.

\subsubsection{Initial Boundary Value Problem}
We will consider Equations~\eqref{eq:2dle1} and~\eqref{eq:2dle2} on a bounded spatial domain $\Omega$, for a time interval $[0,T]$.
For these equations to have a unique solution, initial and boundary conditions must be prescribed.
The resulting initial boundary value problem (IBVP) reads
\begin{align}
    \vspace{-10.0em}
    \boldsymbol{\upsigma}&=\mathbf{C} \boldsymbol{\upvarepsilon} (\boldsymbol{u}) &&\forall \boldsymbol{x} \in \Omega
    \label{eq:leibvp1}
    \\
    \mathbf{D}^{\textsc{t}} \boldsymbol{\upsigma} - \ddot{\boldsymbol{u}} &= - \boldsymbol{f} &&\forall \boldsymbol{x} \in \Omega
    \label{eq:leibvp2}
    \\
    \boldsymbol{u}(t,\boldsymbol{x}) &= \boldsymbol{g}_{\textsc{d}}(t,\boldsymbol{x}) &&\forall \boldsymbol{x} \in{\partial \Omega_{\textsc{d}}}    
    \\
    \boldsymbol{\upsigma}(t,\boldsymbol{x}) \cdot \boldsymbol{n} &= \boldsymbol{g}_{\textsc{n}}(t,\boldsymbol{x}) &&\forall \boldsymbol{x} \in{\partial \Omega_{\textsc{n}}}    
    \\
    \boldsymbol{u}(0,\boldsymbol{x}) &= \boldsymbol{\mu}_0(\boldsymbol{x}) &&\forall \boldsymbol{x} \in \Omega
    \\
    \boldsymbol{\dot{u}}(0,\boldsymbol{x}) &= \boldsymbol{\rho}_0(\boldsymbol{x}) &&\forall \boldsymbol{x} \in \Omega,
\end{align}
where $\partial \Omega_{\textsc{d}}$ and $\partial \Omega_{\textsc{n}}$ are the parts of the boundary with respectively Dirichlet boundary condition
$\boldsymbol{g}_\textsc{d}$ and Neumann boundary condition $\boldsymbol{g}_\textsc{n}$, $\boldsymbol{n}$ is the unit vector normal to the boundary, and $\boldsymbol{\mu}_0$ and $\boldsymbol{\rho}_0$ are the initial value conditions. In the present work, we use Dirichlet boundary conditions on all of the boundaries, i.e., $\partial \Omega_{\textsc{d}} = \partial \Omega$.

\subsubsection{Variational formulation}
Using the Galerkin approach we can transform the IBVP above to a Variational Problem (VP). To facilitate this transformation we introduce the following function spaces:
\begin{align}
    X &= \boldsymbol{H}^1(\Omega) = [H^1(\Omega)]^{n} \label{eq:Hilbert-Space}\\
    X_{\textsc{d}} &= \{ 
    \boldsymbol{v} \in X: \, \boldsymbol{v}=\boldsymbol{g}_{\textsc{d}} \,\, \mathrm{on} \,\, \partial \Omega_{\textsc{d}} \} \\
    X_0 &= \{ \boldsymbol{v} \in X: \, \boldsymbol{v}=0 \,\, \mathrm{on} \,\, \partial \Omega_{\textsc{d}} \}
    \\
    Y(X) &= \{ \boldsymbol{v} : \forall t \in[0,T], \boldsymbol{v}(\boldsymbol{x},t)  \in X, \, \int_0^T ||\boldsymbol{v}||^2_{ \boldsymbol{H}^2(\Omega)} \, dt < \infty \}
    \label{eq:Time-Spatial-Space}
\end{align}
Here, $\boldsymbol{H}^1(\Omega) = [H^1(\Omega)]^{n}$ and $\boldsymbol{H}^2(\Omega) = [H^2(\Omega)]^{n}$ are the Hilbert spaces for functions with first and second order weak derivatives, respectively,  for problems with $\Omega \in \mathbb{R}^{n}$, and we define the spaces $Y(X_{\textsc{d}})$ and $Y(X_0)$ by substituting $X$ with $X_{\textsc{d}}$ and $X_0$ in Equation~\eqref{eq:Time-Spatial-Space}.
Furthermore, we split the unknown displacement into two parts: $\boldsymbol{u} = \boldsymbol{u}_0 + \boldsymbol{u}_{\textsc{d}}$, where $\boldsymbol{u},\boldsymbol{u}_{\textsc{d}} \in X_{\textsc{d}}$ and $\boldsymbol{u}_0 \in X_0$. The so-called {\em lifting\/} function $\boldsymbol{u}_{\textsc{d}}$ is introduced in order to handle non-homogeneous Dirichlet conditions, see e.g.,~\citep{quarteroni2017nummod} for details. The variational formulations then read
\begin{equation}
   \mathrm{Find} \,\, \boldsymbol{u}_0 \in Y(X_0) : \,\, 
   a(\boldsymbol{u}_0,\boldsymbol{v}) + \frac{d^2}{dt^2}(\boldsymbol{u}_0,\boldsymbol{v}) = l(\boldsymbol{v}) 
   \,\,\,\, \forall \boldsymbol{v} \in X_0, 
   \label{eq:VP}
\end{equation}

where the bilinear forms are defined as
\begin{align}
  a(\boldsymbol{w},\boldsymbol{v}) &= \int_{\Omega} \boldsymbol{\upvarepsilon}(\boldsymbol{v})^{\textsc{t}}
                         \mathbf{C}
                         \boldsymbol{\upvarepsilon}(\boldsymbol{w}) \, d \Omega
  \,\,\,\, \forall \boldsymbol{w}, \boldsymbol{v} \in X,
  \label{a-form}
  \\
  (\boldsymbol{w},\boldsymbol{v}) &= \int_{\Omega} \boldsymbol{v}^{\textsc{t}} \boldsymbol{w}  \, d \Omega
  \,\,\,\, \forall \boldsymbol{w}, \boldsymbol{v} \in X,
  \label{m-form}
\end{align}
and the linear form reads
\begin{equation}
   l(\boldsymbol{v}) = \int_\Omega \boldsymbol{f} \boldsymbol{v} \, d \Omega
                     + \int_{\partial \Omega_{\textsc{n}}} \boldsymbol{g}_{\textsc{n}} \boldsymbol{v} \, d \partial \Omega
                     - a(\boldsymbol{u}_{\textsc{d}},\boldsymbol{v})
                     - \frac{d^2}{dt^2}(\boldsymbol{u}_{\textsc{d}},\boldsymbol{v}) 
                    \,\,\,\, \forall \boldsymbol{v} \in X_0,
\label{l-form}                
\end{equation}

\subsubsection{Semi-discretization with Finite Elements}
By discretizing in space we obtain the semi-discrete formulation:
\begin{equation}
   \mathrm{Find} \,\, \boldsymbol{u}_0 \in Y(X_{h,0}) : \,\, 
   a(\boldsymbol{u}_{h,0},\boldsymbol{v}) + \frac{d^2}{dt^2}(\boldsymbol{u}_{h,0},\boldsymbol{v}) = l(\boldsymbol{v}) 
   \,\,\,\, \forall \boldsymbol{v} \in X_{h,0}, 
   \label{eq:FEM-Semidiscrete}
\end{equation}
Here, $X_{h,0} \subset X_0$, and $\mathrm{dim}(X_{h,0}) = N_{\mathrm{dof}} < \infty$.
For linear finite elements we can express $X_{h,0}$ as
\begin{align}
    X_{h,0} &= \left\{ v\in X_0 \, | \, v_{|\Omega_{\scriptstyle e}} \in P_1(\Omega_e), \,\, e=1, \ldots, N_{\mathrm{el}} \right\} \\
            &= \mathrm{span} \{ \phi_1, \phi_2, \ldots, \phi_{N_{\scriptstyle \mathrm{dof}}} \}
\end{align}
where the domain $\Omega$ are properly triangulated into simplex elements (lines in 1D, triangles in 2D, and tetrahedrons in 3D) with domain $\Omega_e$, and $N_{\mathrm{dof}}$ and $N_{\mathrm{el}}$ are the number of degrees of freedom and number of finite elements, respectively. Using the nodal basis functions $\{\phi_i\}, \,\, i=1,\ldots, N_{\mathrm{dof}}$ we can represent any test function in $X_{h,0}$ and trial function in $Y(X_{h,0})$ as follows:
\begin{align}
    \forall \boldsymbol{v} \in X_{h,0} \,\,\,\,\,\,\,\, 
    \boldsymbol{v}(\boldsymbol{x}) &= \sum_{i=1}^{N_{\mathrm{dof}}} v_i \, \phi_i(\boldsymbol{x}) \\
    \forall \boldsymbol{u}_{h,0}(\boldsymbol{x}) \in Y(X_{h,0}) \,\,\,\,\,\,\,\, 
    \boldsymbol{u}_{h,0} &= \sum_{j=1}^{N_{\mathrm{dof}}} (u_{h,0})_j(t) \, \phi_i(\boldsymbol{x})
\end{align}
Notice that, the coefficients for the trial functions, $(u_{h,0})_j(t)$, are time-dependent, whereas that is not the case for the coefficients for the test functions, $v_i$.

By insertion of the nodal basis into Equation~\eqref{eq:FEM-Semidiscrete} we get the following system of ordinary differential equations (ODEs):
\begin{equation}
    \boldsymbol{A} \boldsymbol{U}(t) +  \boldsymbol{M} \ddot{\boldsymbol{U}}(t) = \boldsymbol{F}(t)
    \label{eq:ODE-System}
\end{equation}
where an element of the system {\em stiffness matrix\/} $\boldsymbol{A}$ is defined by $A_{ij}= a(\phi_i,\phi_j)$, an element of the system {\em mass matrix\/} $\boldsymbol{M}$ is defined by $M_{ij}= (\phi_i,\phi_j)$, an element of the system {\em load vector\/}, $\boldsymbol{F}$ is defined by $F_{i}= l(\phi_i)$, and an element of the unknown finite element solution $\boldsymbol{U}$ is $U_i(t) = (u_{h,0})_i(t)$, where the range of the indexes are: $i,j=1,\ldots,N_{\mathrm{dof}}$.

\subsubsection{System of discrete equations}
To obtain a fully discrete system of equations we need to discretize Equation~\eqref{eq:ODE-System} in time. We will here use a second order accurate implicit Euler finite-difference approximations:
\begin{equation}
   \boldsymbol{A} \boldsymbol{U}^{(i+1)} 
   +  \boldsymbol{M} \frac{1}{k^2} \left[ \boldsymbol{U}^{(i-1)} - 2  \boldsymbol{U}^{(i)} +  \boldsymbol{U}^{(i+1)} \right]
   = \boldsymbol{F}^{(i+1)}  
\end{equation}
Here, $\boldsymbol{U}^{(i)}$ denotes the approximate solution $\boldsymbol{U}(t)$ at the time $t^i = ik$, where $k$ is the time step. After rearranging the term we get the following system of algebraic equations:
\begin{equation}
   \left( \boldsymbol{A} + \frac{1}{k^2}\boldsymbol{M} \right) \boldsymbol{U}^{(i+1)} 
   = \boldsymbol{F}^{(i+1)} + \frac{1}{k^2} \boldsymbol{M} \left[ 2  \boldsymbol{U}^{(i)} -  \boldsymbol{U}^{(i-1)} \right]
\label{eq:Algebraic-System}     
\end{equation}
This system can be solved by any appropriate solver, preferably a sparse solver for 2D and small 3D problems. For large problems (e.g., in 3D) iterative solvers as Conjugate Gradient Method will typically be a better choice. Starting from $U_0$, the known initial condition $\boldsymbol{\mu}_0$ projected onto $X_h$, repeatedly solving Equation~\eqref{eq:Algebraic-System} allows us to march forward in time, eventually obtaining approximations of $\boldsymbol{u}(\boldsymbol{x},t)$ up to the specified final time $T$. For the spatial discretization, we use piecewise linear Lagrange simplex elements~\citep{Courant1943VariationalMF} on an equidistant grid. As can be seen in Section~\ref{sec:res:el}, the relative error of the resulting model's predictions is roughly 1\% or smaller in the scenarios where no modeling error is synthesized. This is a reasonable result for what must be considered as a fairly simple PBM, and serves as a good benchmark for DDM and CoSTA to beat.

\subsection{Nonlinear Elasticity Modeling}
\label{sec:nonlin_elasticity}
As previously mentioned, the linear elasticity model introduced in the previous section is well-suited for modeling small strains and stresses.\footnote{The precise meaning of ``small'' strains and stresses is highly material-dependent.} More flexible elasticity modeling can be achieved by relaxing the assumptions that material properties like Young's modulus $E$ and the Poisson ratio $\nu$ be constants. For the nonlinear elasticity cases considered herein, we consider $E$ as a function of strain, i.e.\ $E = E(\boldsymbol{\upvarepsilon})$. To accommodate for this, a slight alteration to the system \eqref{eq:2dle1}--\eqref{eq:2dle2} is required, such that we obtain
\begin{align}
    \boldsymbol{\upsigma}&=\mathbf{C}(\boldsymbol{\upvarepsilon}) \boldsymbol{\upvarepsilon}
    \label{eq:2dn1e} \\
    \mathbf{D}^{\textsc{t}} \boldsymbol{\upsigma} - \ddot{\boldsymbol{u}} &= - \boldsymbol{f}.
    \label{eq:2dn2e}
\end{align}
Here, we highlight that Equation~\eqref{eq:2dn1e} is nonlinear in $\boldsymbol{u}$ via the Young's modulus $E$, and hence the constitutive matrix $\mathbf{C}$,  dependence on the strain vector $\boldsymbol{\varepsilon}$. 

\subsection{Data Driven modeling Using Neural Networks}
\label{subsec:ddm}

We use deep neural networks (DNN) to create the purely data-driven models considered in herein. This is motivated by the observation that DNNs are universal approximators~\citep{Cybenko1989ApproximationBS, Hornik1989MultilayerFN, LESHNO1993861} and hence have the ability to model highly complex nonlinear phenomena. 
A DNN with $N$ 
layers is a function $D: \mathbb{R}^{d_0} \rightarrow \mathbb{R}^{d_N}$ defined by
\begin{equation}
D(x) = \mathcal{T}_N \circ s_{N-1} \circ \mathcal{T}_{N-1} \circ s_{N-2} \cdots \circ s_1 \circ \mathcal{T}_1 (x)
\label{eq:DNN}
\end{equation}
where $\mathcal{T}_i: \mathbb{R}^{d_{i-1}} \rightarrow \mathbb{R}^{d_i}$ are affine transformations and $s_i: \mathbb{R}^{d_{i}} \rightarrow \mathbb{R}^{d_i}$ are some preferably nonlinear activation functions.\footnote{With linear activation functions the method reduces to multivariate linear regression. For the universal approximation property it should also be nonpolynomial.} Each transformation $\mathcal{T}_i$ is determined by a weight matrix and a bias vector, a total of $d_{i-1} \cdot d_{i} + d_i$ values. These values are tuned to minimize the DNN's predictive error on the training data. DNNs are widely used, due to their simplicity and yet astonishing performance in many situations. They can be applied to a large variety of problems, and have delivered impressive achievements across numerous fields of research and applications. Meanwhile, they also have some inherent weaknesses. They are prone to overfitting on the training data. While there are many ways of preventing the network from being too specialized~\citep{sym10110648}, they will not be good at extrapolation cases where the prediction task is somehow qualitatively different from the training tasks. In addition, they are not easily explainable or predictable, meaning it's hard to explain what kind of patterns the network will look for, and unexpected predictions can occur. Compared to PBMs, DDMs usually make predictions much faster, but may need long training times for optimal performance.

\subsection{Corrective Source Term Approach}
    \label{subsec:costa}
    In this section we present a brief justification of the use of CoSTA. It is based on the more elaborate argument that can be found in~\cite{Blakseth2022dnn}.
    
    Consider the differential equations 
    \begin{align}
        L_{\Omega} \omega&=\theta, \quad \forall \boldsymbol{x} \in \Omega \nonumber \\ 
        L_{\partial \Omega}\omega&=\psi, \quad \forall \boldsymbol{x} \in \partial \Omega,
        \label{eq:PDE}
    \end{align}
    where $L_{\Omega}, L_{\partial \Omega}$ are differential operators, $\theta$ is a source term, $\psi$
    a function specifying the boundary condition, and $\omega$ is the unknown of interest. For notational simplicity we assume $\omega$ to be a scalar. 
    Now let $\tilde{\omega}$ be the solution to the perturbed problem
    \begin{align}
        \tilde{L}_{\Omega}\tilde{\omega}&=\tilde{\theta}, \quad \forall\boldsymbol{x} \in \Omega \nonumber\\
        \tilde{L}_{\partial \Omega} \tilde{\omega}&=\tilde{\psi}, \quad \forall \boldsymbol{x} \in \partial \Omega,
        \label{eq:perturbed_PDE}
    \end{align}
    where the perturbations $\tilde{\cdot}$ are due to imperfections such as unknown physics, modeling errors, discretization error, or noise. For example, we will later in this paper simplify a nonlinear operator $L_\Omega$ with a linear (and discretized) $\tilde{L}_{\Omega}$.
    Assume we can calculate the residuals defined as
    \begin{align}
        r_{\Omega} &= \tilde{L}_{\Omega}(\omega-\tilde{\omega}) \nonumber\\
        r_{\partial \Omega} &=\tilde{L}_{\partial \Omega}(\omega-\tilde{\omega}),
    \end{align}
    and let $\tilde{\tilde{\omega}}$ be the solution of the corrected, perturbed problem
    \begin{align}
        \tilde{L}_{\Omega}\tilde{\tilde{\omega}}&=\tilde{\theta}+r_{\Omega}, \quad \forall \boldsymbol{x} \in \Omega \nonumber
        \\
        \tilde{L}_{\partial \Omega} \tilde{\tilde{\omega}}&=\tilde{\psi}+r_{\partial \Omega}  \quad \forall \boldsymbol{x} \in \partial \Omega.
        \label{eq:corrected_PDE}
    \end{align}
    Using the definition of the residuals, as well as the perturbed Equation~\eqref{eq:perturbed_PDE}, we see that 
    \begin{align*}
        \tilde{L}_{\Omega}\tilde{\tilde{\omega}}&=\tilde{\theta}+\tilde{L}_{\Omega}(\omega-\tilde{\omega}) = \tilde{L}_{\Omega}\omega, \quad \forall \boldsymbol{x} \in \Omega\\
        \tilde{L}_{\partial \Omega} \tilde{\tilde{\omega}}&=\tilde{\psi}+\tilde{L}_{\partial \Omega}(\omega-\tilde{\omega}) = \tilde{L}_{\partial \Omega}\omega  \quad \forall \boldsymbol{x} \in \partial \Omega,
    \end{align*}
    which reduces to $\tilde{\tilde{\omega}}=\omega$ if the corrected, perturbed problem~\eqref{eq:corrected_PDE} 
    yields a unique solution.
    From this argument, we see that the source term corrections are able to compensate for perturbations in the differential operators as well as the source terms.

    In real scenarios, we obviously cannot calculate the residual exactly, as that would require the solution we are trying to estimate. The idea of the CoSTA method is to use a DDM to estimate the residual. For the input of the DDM we use the uncorrected solution $\tilde{\omega}$. This means the PBM is used twice per time level, first to solve Equation~\eqref{eq:perturbed_PDE} for $\tilde{\omega}$, then Equation~\eqref{eq:corrected_PDE} for $\tilde{\tilde{\omega}}$. As mentioned in \cite{Blakseth2022dnn}, other DDM input choices are possible. For example, one could use the corrected prediction $\tilde{\tilde{\omega}}$ from the previous time level, such that the PBM is only used once per time level. Investigating different DDM input choices is outside the scope of the present work, so we stick with the DDM input used in \cite{Blakseth2022dnn}.

    In the case of linear elasticity, the variational formulation of Equation~\eqref{eq:perturbed_PDE} corresponds to Equation~\eqref{eq:VP}. The corrected form of Equation~\eqref{eq:Algebraic-System}, corresponding to the general corrected perturbed Equation~\eqref{eq:corrected_PDE}, reads
    \begin{equation}
   \left( \boldsymbol{A} + \frac{1}{k^2}\boldsymbol{M} \right) \boldsymbol{U}^{(i+1)} 
   = \boldsymbol{F}^{(i+1)} + \frac{1}{k^2} \boldsymbol{M} \left[ 2  \boldsymbol{U}^{(i)} -  \boldsymbol{U}^{(i-1)} \right] + \boldsymbol{r}
   \label{eq:CoSTA_elasticity}     
\end{equation}
    where $\boldsymbol{r}$ is the residual vector corresponding to Equation~\eqref{eq:Algebraic-System} with an appropriate projection of the true solution inserted in the place of $U$.
    In the present work, we use a DNN to approximate the residual vector $\boldsymbol{r}$.
    As the boundary values are known (as we here assume $\partial \Omega_{\textsc{d}} = \partial \Omega$), we do not need any correction for these elements. Therefore the DNN only needs to output values for the interior nodes.

\subsection{Method of Manufactured Solutions}
    \label{subsec:manufacturedsolution}
    We seek to evaluate the performance of the models described above on a selection of different elasticity modeling problems. To this end, we use the method of manufactured solutions~\citep{Roache2001cvb} to create exact reference data. For a governing equation written on the general form~\eqref{eq:PDE}, the method involves choosing a solution $\omega$, and calculating the source $\theta$ that admits the chosen $\omega$ as a solution to the governing equation. In order to evaluate the accuracy of a model, the model is used to approximate the chosen solution $\omega$ given the calculated $\theta$ as well as the correct boundary data $\omega|_{t=0}$ and $\omega|_{\Omega}$.
    
    The alternative to using the method of manufactured solutions is to choose the conditions $\theta$, $\omega|_{t=0}$ and $\omega|_{\Omega}$ to use for approximating the solution $\tilde{\omega}$. But then the correct solution is unknown, so the error $\omega-\tilde{\omega}$ must be approximated by using a more precise method.\footnote{usually high fidelity numerical solutions of the equation, with a fine grid.} Usually, this is very computationally demanding, and it only yields an approximate error.
    In comparison, the method of manufactured solutions is computationally inexpensive and introduces no uncertainty in the assessment of model accuracy.

\section{Methodology}
\label{sec:methodology}
In the present work, we consider four different numerical experiments, simply enumerated as Experiment 1--4. The goal of these experiments is to test the predictive accuracy and generalizability of CoSTA, in comparison to its constituent PBM and DDM components. This section is dedicated to describing the methodology of the experiments. Elements of the methodology that are common among all experiments are described in Section~\ref{sec:genset}. What is unique for each experiments is then described in Sections~\ref{sec:el}--\ref{sec:elnl}. This is summarized in Table~\ref{tab:experiments}. All four experiments concern elasticity modeling using the PBM, DDM and CoSTA models introduced earlier. 
\begin{table*}
    \centering
        \caption{Description of experiments and modeling errors (in addition to discretization error).}
    \label{tab:experiments}
    \begin{tabular}{ccc}
        \toprule
        Exp.\ \# & Modelled physics & PBM modeling error \\
        \hline
        1 & 2D Linear elasticity with load term & None \\
        2 & 2D Linear elasticity with load term & Load term replaced with zero \\
        3 & 3D Linear elasticity with reduced dimensionality & A dimension is ignored \\
        4 & 2D Nonlinear elasticity & PDE is linearized \\
        \bottomrule
    \end{tabular}
\end{table*}

\subsection{General Setup}
\label{sec:genset}

\subsubsection{Data Generation}
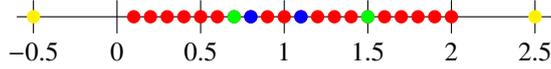
\begin{figure}
    \centering
  \begin{tikzpicture}[scale=2.2]
    \draw (-0.6,0) -- (2.6,0);
    \foreach \i in {-0.5,0,...,2.5} 
      \draw (\i,0.1) -- + (0,-0.2) node[below] {$\i$}; 
    \foreach \i in {0.1, 0.2,..., 2.0}
      \fill[red] (\i,0) circle (0.4 mm);
    \foreach \i in {0.8, 1.1}
      \fill[blue] (\i,0) circle (0.4 mm);
    \foreach \i in {0.7, 1.5}
      \fill[green] (\i,0) circle (0.4 mm);
    \foreach \i in {-0.5, 2.5}
      \fill[yellow] (\i,0) circle (0.4 mm);
  \end{tikzpicture}
    \caption{
    Values of $\alpha$ used for \disc{red}training, \disc{blue}validation, \disc{green}interpolation testing and \disc{yellow}extrapolation testing.
    }
    \label{fig:alphas}
\end{figure}

To conduct our numerical experiments, we need data for two purposes. First of all, we need reference data to which we can compare the predictions of our models. Moreover, we need training data for the purely data-driven model and the data-driven component of CoSTA. To obtain this data, we use the method of manufactured solutions (cf.\ Section~\ref{subsec:manufacturedsolution}).\footnote{In a real-world use-case, one might instead use sensor measurements or simulation data from a high-fidelity offline model.} In order to cover a broad class of solutions, the manufactured solutions include a variety of polynomial, exponential and harmonic functions.
These functions are described in the sections covering experiment-specific information. What they all have in common is that they are parametrized by some variable $\alpha$. Consequently, it is useful to view each manufactured solution as a family of functions, one for each value of $\alpha$, corresponding to different system states. In the present work, we consider a total of 22 $\alpha$-values, as listed in Table~\ref{table:alphas} and visualized in Figure~\ref{fig:alphas}. 16 of these values are used for training the DNNs used in the DDM and CoSTA models,\footnote{We highlight that PBM does not require any training.} while two more are used for DNN validation. The remaining four values are used for testing the accuracy of the models. Two of these lie within the range used to generate training data, while the other values are outside this range. We refer to the former as interpolation scenarios and the latter as extrapolation scenarios.

As mentioned in Section~\ref{subsec:lin_elasticity}, we have two quantities of main interest in elasticity modeling: $\boldsymbol{\upsigma}$ and $\boldsymbol{u}$.\footnote{Recall that $\boldsymbol{\sigma}$ and $\boldsymbol{\upsigma}$ are equivalent representations of the same variables.} These are related through the constitutive relation~\eqref{eq:2dle1}, which depends on the material properties $E$ and $\nu$. Unless otherwise stated, we always use $E=1$ and $\nu=0.25$. Using the constitutive relation, if we know $\boldsymbol{u}$, we can compute $\boldsymbol{\upsigma}$, and vice versa. For our manufactured solutions, we elect to fix $\boldsymbol{u}$, $E$ and $\nu$. $\boldsymbol{\sigma}$ is then computed using Equation~\eqref{eq:2dle1}, before $\boldsymbol{f}$ is computed using Equation~\eqref{eq:2dle2}.

\begin{table}
    \centering
    \caption{Values of the parameter $\alpha$ used for training, validation and testing.}
    \begin{tabular}{ ccc } 
    \toprule
    Set usage &Notation& Values \\
    \hline
    Testing &$\mathcal{A}_{\text{test}}$& $\{-0.5,0.7,1.5,2.5\}$ \\ 
    Validation &$\mathcal{A}_{\text{val}}$& $\{0.8,1.1\}$ \\ 
    Training &$\mathcal{A}_{\text{train}}$& $\{0.1,0.2,... 2.0\} \setminus (\mathcal{A}_{\text{test}} \cup \mathcal{A}_{\text{val}})$ \\
    \bottomrule
    \end{tabular}
    \label{table:alphas}
\end{table}

\subsubsection{Model Summary}
In our numerical experiments, we use the PBM, DDM and CoSTA models for elasticity problems that were introduced in Section~\ref{sec:theory}. These are briefly summarized below.

\textbf{PBM:} Equation~\eqref{eq:Algebraic-System} is solved using a linear finite element model with first order (i.e.\ piecewise linear) simplex Lagrange elements~\citep{Courant1943VariationalMF} on an equidistant grid. Unless otherwise specified, the spatial domain considered is the unit square $[0,1]\times [0,1]$, and the time domain is the unit interval $[0,1]$. The time interval is divided into time steps of constant length $k=K^{-1}$, where $K$ is the number of time steps.
The spatial and temporal resolution used varies across experiments, as shown in Table~\ref{tab:resolutions}.
\begin{table}
    \centering
        \caption{Overview of spatial and temporal resolution.}
    \label{tab:resolutions}
    \begin{tabular}{ccc}
        \toprule
        Exp.\ \# & No.\ Elements & No. Time Steps, $K$ \\
        \hline
        1 & $15 \times 15$ & $1000$\\
        2 & $15 \times 15$ & $1000$\\
        3 & $15 \times 15$ & $1000$\\
        4 & $10 \times 10$ & $500$\\
        \bottomrule
    \end{tabular}
\end{table}

\textbf{DDM:} For the DDM, we use a DNN with four dense hidden layers of 80 nodes each. The length of the input and output layers depend on the PDE and discretization. The input vector contains every basis function, while the output does not contain the functions on the boundary.
Recall that the two-dimensional (2D) elasticity equation has a vector field solution with two basis functions per element.	
For the experiments with 15 elements in each dimension, this gives $(2\cdot16)^2$ input nodes and $(2\cdot14)^2$ output nodes, as shown in Figure~\ref{fig:arch}.

\begin{figure*}
    \centering
    \def\layersep{2.2cm}
    \begin{tikzpicture}[shorten >=1pt,->,draw=black!50, node distance=\layersep]
        \tikzstyle{every pin edge}=[<-,shorten <=1pt]
        \tikzstyle{neuron}=[circle,fill=black!25,minimum size=17pt,inner sep=0pt]
        \tikzstyle{input neuron}=[neuron, fill=green!50];
        \tikzstyle{input2 neuron}=[neuron, fill=yellow!50];
        \tikzstyle{output neuron}=[neuron, fill=red!50];
        \tikzstyle{hidden neuron}=[neuron, fill=blue!50];
        \tikzstyle{annot} = [text width=4em, text centered]
        \foreach \name / \y in {1,2,3,5}
            \node[input neuron, pin=left:Input] (I-\name) at (0,-\y) {};
        \node at (0,-4 cm) {$\vdots$};
        \foreach \name / \y in {1,2,3,4,6}
            \path[yshift=0.5cm]
                node[hidden neuron] (H1-\name) at (\layersep,-\y cm) {};
        \path[yshift=0.5cm]node at (\layersep,-5 cm) {$\vdots$};
        \foreach \name / \y in {1,2,3,4,6}
            \path[yshift=0.5cm]
                    node[hidden neuron] (H2-\name) at (2*\layersep,-\y cm) {};
        \path[yshift=0.5cm]node at (2*\layersep,-5 cm) {$\vdots$};
        \foreach \name / \y in {1,2,3,4,6}
            \path[yshift=0.5cm]
                node[hidden neuron] (H3-\name) at (3*\layersep,-\y cm) {};
        \path[yshift=0.5cm]node at (3*\layersep,-5 cm) {$\vdots$};
        \foreach \name / \y in {1,2,3,4,6}
            \path[yshift=0.5cm]
                node[hidden neuron] (H4-\name) at (4*\layersep,-\y cm) {};
        \path[yshift=0.5cm]node at (4*\layersep,-5 cm) {$\vdots$};
        \foreach \name / \y in {1,2,4}
            \path[yshift=-0.5cm]
                node[output neuron,pin={[pin edge={->}]right:Output}, right of=H3-3] (O-\name) at (4*\layersep,-\y) {};
        \path[yshift=-0.5cm]
            node[right of=H3-3] (O-3) at (4*\layersep,-3) {$\vdots$};
        \foreach \source in {1,2,3,5}
            \foreach \dest in {1,2,3,4,6}
                \path (I-\source) edge (H1-\dest);
        \foreach \source in {1,2,3,4,6}
            \foreach \dest in {1,2,3,4,6}
                \path (H1-\source) edge (H2-\dest);
        \foreach \source in {1,2,3,4,6}
            \foreach \dest in {1,2,3,4,6}
                \path (H2-\source) edge (H3-\dest);
        \foreach \source in {1,2,3,4,6}
            \foreach \dest in {1,2,3,4,6}
                \path (H3-\source) edge (H4-\dest);
        \foreach \source in {1,2,3,4,6}
            \foreach \dest in {1,2,4}
                \path (H4-\source) edge (O-\dest);
        \node[annot,above of=H1-1, node distance=1cm] (hl1) {Hidden layer};
        \node[annot,right of=hl1] (hl2) {Hidden layer};
        \node[annot,right of=hl2] (hl3) {Hidden layer};
        \node[annot,right of=hl3] (hl4) {Hidden layer};
        \node[annot,left of=hl1] {Input layer};
        \node[annot,right of=hl4] {Output layer};
        \node[annot,below of=H1-6, node distance=1cm] (n1) {80 nodes};
        \node[annot,right of=n1] (n2) {80 nodes};
        \node[annot,right of=n2] (n3) {80 nodes};
        \node[annot,right of=n3] (n4) {80 nodes};
	    \node[annot,left of=n1] {$(2\cdot16)^2$ nodes};	
        \node[annot,right of=n4] {$(2\cdot14)^2$ nodes};	
    \end{tikzpicture}	
    \caption{Visualization of the DNN architecture for the experiments with $15\times 15$ elements. 
    The nodes represent input, output and intermediate values, while the arrows going between them represent dependencies.
    Generally there is one input node for each basis function, and one output node for each basis function not on the edge.
    }
    \label{fig:arch}
\end{figure*}

The neural networks are implemented using TensorFlow~\citep{tensorflow2015-whitepaper}. As activation function we use leaky~ReLU~\citep{lrelu}, with coefficient 0.01 for negative inputs. Training parameters are presented in Table~\ref{table:params}.
A patience of 20 means the training stops when the score on the \textit{validation} set has not improved for the last 20 optimizer steps.
All the data is normalized before it is inputted in the DNN, and the output is unnormalized, based on the statistical properties of the training data.


\textbf{CoSTA:} Equation~\eqref{eq:CoSTA_elasticity} is solved using the exact same finite element method discretization used in the PBM with exactly the same spatial and temporal resolution. Moreover, the residual $r$ is approximated by a DNN with exactly the same architecture and hyperparameter values as the DNN used for pure DDM.

\subsubsection{Experimental Procedure}

The experimental procedure used is presented in Algorithm~\ref{alg:dimred}. $K$ denotes the total number of time steps.
Note the difference between $\check{\boldsymbol{U}}$, $\bar{\boldsymbol{U}}$ and $\hat{\boldsymbol{U}}$: we denote by $\check{\boldsymbol{U}}$ the vector form of the exact solution\footnote{That is, the piecewise linear function characterised by $\check{\boldsymbol{U}}$ equals the exact solution on the grid points}. The predictions $\bar{\boldsymbol{U}}$, used for training, are based on the previous \textit{exact} step, while $\hat{\boldsymbol{U}}$, used for testing, are based on the previous \textit{predicted} step (the first predicted step is based on the exact initial values).
As in Section~\ref{sec:theory}, we use the superscript to denote the time step, e.g.\ $\hat{\boldsymbol{U}}^{(i)}$ is at time $t_i = ik$.

\begin{algorithm}
\caption{Pseudocode showing how the experiments were performed.}
\label{alg:dimred}
    Pick a solution $\boldsymbol{u}_{\mathrm{exact}}(t,x,y,z,\alpha)$\;
    Use governing equations (the PDE, before any simplification) to calculate $\boldsymbol{f}$ from $\boldsymbol{u}_{\mathrm{exact}}$\;
    \For{$\alpha \in \mathcal{A}_{\text{train}}$ and $\alpha \in \mathcal{A}_{\text{val}}$}
    {
        \For{$i=0,1,2...K-1$}{Use (simplified) PBM to calculate $\bar{\boldsymbol{U}}_\text{PBM}^{(i+1)} = \text{PBM}(\check{\boldsymbol{U}}^{(i)})$\;
        Calculate the residual $\boldsymbol{r}^{(i+1)} = \hat{L}(\check{\boldsymbol{U}}^{(i+1)} - \bar{\boldsymbol{U}}_{\text{PBM}}^{(i+1)})$;}
    }
    Train DDM to map $\check{\boldsymbol{U}}^{(i)}$ to $\check{\boldsymbol{U}}^{(i+1)}$, i.e. minimize 
    \\ \hspace*{10mm}$\sum_{\alpha \in \mathcal{A}_{\text{train}}}\sum_{i=0}^{K-1}\left|\check{\boldsymbol{U}}^{(i+1)} - \text{DDM}(\check{\boldsymbol{U}}^{(i)})\right|^2$\;
    Train CoSTA network to map $\bar{\boldsymbol{U}}_{\text{PBM}}^{(i+1)}$ to $\boldsymbol{r}^{(i+1)}$, i.e. minimize 
    \\ \hspace*{10mm}$\sum_{\alpha \in \mathcal{A}_{\text{train}}}\sum_{i=0}^{K-1}\left| \boldsymbol{r}^{(i+1)}-\text{DDM}_{\text{CoSTA}}(\text{PBM}(\check{\boldsymbol{U}}^{(i)}))
    \right|^2$\;
    $\hat{\boldsymbol{U}}_\text{DDM}^{(0)} = \check{\boldsymbol{U}}^{(0)}$\;
    $\hat{\boldsymbol{U}}_\text{PBM}^{(0)} = \check{\boldsymbol{U}}^{(0)}$\;
    $\hat{\boldsymbol{U}}_\text{CoSTA}^{(0)} = \check{\boldsymbol{U}}^{(0)}$\;
    \For{$\alpha \in \mathcal{A}_{\text{test}}$}
    {
        \For{$i=0,1,2...K-1$}{calculate $\hat{\boldsymbol{U}}_\text{DDM}^{(i+1)} = \text{DDM}(\hat{\boldsymbol{U}}_\text{DDM}^{(i)})$\;
        calculate $\hat{\boldsymbol{U}}_\text{PBM}^{({i+1})} = \text{PBM}(\hat{\boldsymbol{U}}_\text{PBM}^{(i)})$\;
        calculate $\hat{\boldsymbol{U}}_\text{CoSTA}^{(i+1)} = \text{PBM}(\hat{\boldsymbol{U}}_\text{CoSTA}^{(i)}, \text{DDM}_{\text{CoSTA}}(\text{PBM}(\hat{\boldsymbol{U}}_\text{CoSTA}^{(i)})))$;}
    }
    Evaluate by comparing to exact solution $\boldsymbol{u}_{\mathrm{exact}}$
\end{algorithm}

The PBM, DDM and CoSTA models map their predictions $\hat{\boldsymbol{U}}^{(i-1)}$
to $\hat{\boldsymbol{U}}^{(i)}$. This map is used iteratively to calculate $\hat{\boldsymbol{U}}^{(K)}$
at the final time step, from the initial (exact) input $\check{\boldsymbol{U}}^{(0)}$. We should therefore expect that any error in one step propagates into the next. It is not feasible to train DDM or CoSTA to counteract this by fixing the global error, because 1) The models would depend on the arbitrary starting time, and 2) the required amount of data and training time would increase drastically since the number of training examples per time series would decrease from $K$ to 1.
Therefore, the exact solution $\check{\boldsymbol{U}}$ is used as both the input and output training data for every time step.

For testing scenarios, the exact solution at time $t=0$ is used as the input for the first step, and then each model predict each step based on the previous. Boundary conditions are also used as input for each step, along with the source function values for the PBM. The error is measured at each step, and presented
in the result sections.

For actual usage of this method for prediction, the first two lines and the last line of Algorithm~\ref{alg:dimred} is skipped. Instead, $\boldsymbol{f}$ must be known (or approximated), along with $\boldsymbol{u}_\mathrm{exact}$ in training scenarios, and the initial conditions and boundary conditions in the testing scenarios.
\begin{table}
    \centering
    \caption{Parameters used for the training procedures for all of the neural networks}
    \begin{tabular}{ cc } 
    \toprule
    Loss function & MSE\\
    Optimizer & Adam~\citep{kingma2014aam} \\ 
    Learning rate & $1e-5$\\ 
    Patience &  20\\
    \bottomrule
    \end{tabular}
    \label{table:params}
\end{table}
\subsubsection{Evaluation and Visualization}
\label{subsec:eval}

To evaluate the performance of the different methods, relative root mean square error (RRMSE) is measured at each time step. This value is defined
\begin{equation}
    \text{RRMSE}(\boldsymbol{U}^{(i)}, \check{\boldsymbol{U}}^{(i)}) = \frac{||\boldsymbol{U}^{(i)}-\check{\boldsymbol{U}}^{(i)}||_2}{||\check{\boldsymbol{U}}^{(i)}||_2}
\end{equation}
for a prediction $\boldsymbol{U}$ and correct solution $\check{\boldsymbol{U}}$.
For models with stochastic results (i.e. DDM and CoSTA, that includes a DNN with random initialisation), 10 models are trained and used. 
In a real world application, the objective would determine if we are interested in the least amount of error at all time steps, or only at the last one. In this work we are interested in both.
In the result plots presented in Section~\ref{sec:results}, the mean of the RRMSE of the 10 models is plotted as a line on a logarithmic scale, as a function of time steps. The uncertainty is also quantified and visualized in the plots by shading the area between the mean RRMSE, and the mean RRMSE plus one standard deviation.\footnote{The standard deviation is calculated using one reduced degree of freedom due to the estimation of the mean.} Mean RRMSE \textit{minus} one standard deviation is not plotted since it might be negative, which does not work well on logarithmic scales. 

\subsection{Modeling Linear Elasticity With Known and Unknown Load Term}
\label{sec:el}
In Experiments~1 and~2, we consider a system governed by the 2D transient linear elasticity equations~\eqref{eq:2dle1} and~\eqref{eq:2dle2} with a non-zero load term $\boldsymbol{f}$. In Experiment~1, $\boldsymbol{f}$ is known to the PBM and CoSTA models. Consequently, discretization error will be the only source of error in the PBM, and the task of the DNN used in the CoSTA model will be to reduce this error. In Experiment~2, $\boldsymbol{f}$ is assumed unknown, so we set $\boldsymbol{f}=\boldsymbol{0}$ in both the PBM and CoSTA models, thereby synthesizing modeling error.\footnote{In real-world use-cases, one should use the best available \textit{a priori} estimation of $\boldsymbol{f}$.}
Since the DDM does not have explicit knowledge about $\boldsymbol{f}$, there is no difference between Experiments~1 and~2 from the perspective of the DDM.
\begin{table}
    \centering
    \caption{Manufactured solutions for the elastic problems.}
    \begin{tabular}{ cc } 
    \toprule
    Label & $\boldsymbol{u}(t,x,y,\alpha)$ \\ 
    \hline
    \vspace{8pt}
    \textit{e1} & $\left[
    \begin{array}{c}
        \sin(\pi(x+\alpha y))\cos(\alpha t) \\
        \cos(\pi(x+\alpha y))\sin(\alpha t)
    \end{array}
\right]$ \\ 
    \vspace{8pt}
    \textit{e2} & $\left[
    \begin{array}{c}
        \exp(\frac{(-tx^2 + y^2)}{(1+\alpha+t^2)}) \\
        \exp(\frac{(+tx^2 - y^2)}{(1+\alpha+t^2)})
    \end{array}
\right]$ \\ 
    \vspace{4pt}
    \textit{e3} & 
    $\left[
    \begin{array}{c}
        x^3 + y^2(t+0.5)^{1.5} + xy\alpha \\
        x^2 + y^3(t+0.5)^{1.1}+ xy\alpha
    \end{array}
\right]$  \\
    \bottomrule
    \end{tabular}
    \label{table:manufact_sols:2d_elastic}
\end{table}

\subsection{Modeling Linear Elasticity with Reduced Dimensionality}
\label{sec:eldr}

In general, a model's computational complexity and expense increase greatly with the number of dimensions being modelled. Conversely, if the situations allow it, reducing the dimensionality of a model can greatly reduce the model's computational cost. However, the model's accuracy will generally also be reduced. Experiment~3 is designed to investigate whether CoSTA can be used to correct errors introduced by dimensionality reduction. To this end, we use our 2D linear elasticity PBM, DDM and CoSTA models to model the displacement of a 2D plane in a 3D object, as illustrated in Figure~\ref{fig:dimred_cube}. The manufactured solutions, which are solutions of the \emph{3D} linear elasticity equations, are listed in Table~\ref{table:manufact_sols:dimred_elastics}. Note that, since the displacement $\boldsymbol{u}$ is a vector field with as many components as there are dimensions, our 2D models cannot predict the $z$-component of $\boldsymbol{u}$.\footnote{An alternative approach to dimensional reduction for the linear elasticity equation (and other vector field PDEs), is to base the PBM on the three-dimensional model, but replace the derivatives in the ignored z-direction with zero (or some other appropriate value). The resulting PBM would produce predictions for all components of $\boldsymbol{u}$. Although the PBM prediction of the third component likely would be quite inaccurate, the CoSTA term could help. Due to the inter-dependencies of the components, this could potentially give a useful prediction of the third component of the displacement, that in turn could make the other components more accurate. This idea is not pursued further in this paper.} Therefore, we consider only the $x$- and $y$-components of $\boldsymbol{u}$ when evaluating the predictive accuracy of our models. This approach is most relevant when displacement in the reduced direction is either small and/or of little interest compared to displacement in the other directions.
\newcommand{\Depth}{3}
\newcommand{\Height}{3}
\newcommand{\Width}{3/2}
\begin{figure}[htb]
    \centering

    \begin{tikzpicture}
        \coordinate (O) at (0,0,0);
        \coordinate (A) at (0,\Width,0);
        \coordinate (B) at (0,\Width,\Height);
        \coordinate (C) at (0,0,\Height);
        \coordinate (D) at (\Depth,0,0);
        \coordinate (E) at (\Depth,\Width,0);
        \coordinate (F) at (\Depth,\Width,\Height);
        \coordinate (G) at (\Depth,0,\Height);
        \coordinate (H) at (0, -\Width, 0);
        \coordinate (I) at (\Depth, -\Width, 0);
        \coordinate (J) at (\Depth, -\Width, \Height);
        \coordinate (K) at (0, -\Width, \Height);
        \draw[black] (I) -- (J) -- (K);

        \fill[pattern=north west lines,opacity=0.3](O) -- (C) -- (G) -- (D) -- cycle;
        \draw[line width=0.8mm,red] (C) -- (G) -- (D);
        \draw[dashed, line width=0.8mm, color=red] (C) -- (O) -- (D);
        \draw[black] (I) -- (E) -- (A) -- (B) -- (K) -- (B) -- (F) -- (E) -- (F) -- (J);

        \draw[dashed] 
          (A) -- (H) -- (I);
        \draw[dashed]
          (H) -- (K);
        \draw[->] 
          (C) -- ++ (0,0,1) node[anchor=north east]{$x$};
        \draw[->]
           (D) -- ++(1,0,0) node[anchor=north west]{$y$};
        \draw[->]  
          (A) -- ++(0,1,0) node[anchor=south]{$z$};
    \end{tikzpicture}
    \caption{Example of an object to model using the dimensional reduction method. While the full model is a three-dimensional cube, we only predict the shaded plane. Boundary conditions needed are those at the thick red edges at $z=0$.}
    \label{fig:dimred_cube}
\end{figure}
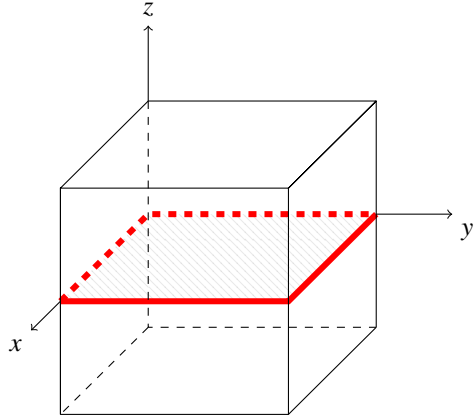

\begin{table}
    \centering
    \caption{Manufactured solutions for the dimensionally reduced elasticity problems.}
    \begin{tabular}{ cc } 
    \toprule
    Label & $u(t,x,y,\alpha)$ \\ 
    \hline
    \vspace{8pt}
    \textit{ed1} & $\left[
    \begin{array}{c}
        \sin(\pi(x+\alpha y + \frac{1+\alpha}{2}z))\cos(\alpha t) \\
        \cos(\pi(x+\alpha y + \frac{1+\alpha}{2}z))\sin(\alpha t) \\
        -\cos(\pi(x+\alpha y + \frac{1+\alpha}{2}z))\sin(\alpha t)
        
    \end{array}
\right]$ \\ 
    \vspace{8pt}
    \textit{ed2} & $\left[
    \begin{array}{c}
        \exp(\frac{(-tx^2 + y^2 + z^2)}{(1+\alpha+t^2)}) \\
        \exp(\frac{(+tx^2 - y^2 + z^2)}{(1+\alpha+t^2)}) \\
        \exp(\frac{(+tx^2 + y^2 - z^2)}{(1+\alpha+t^2)})
    \end{array}
\right]$ \\ 
    \vspace{4pt}
    \textit{ed3} & 
    $\left[
    \begin{array}{c}
        x^3 + y^2(t+\frac{1}{2})^{1.5} + xy\alpha + \sqrt{t+\frac{1}{2}}z^2 + z(x+y)\alpha\\
        x^2 + y^3(t+\frac{1}{2})^{1.1} - xy\alpha + \sqrt{t+\frac{1}{2}}z^2 + z(x-y)\alpha\\
        x^2 + y^2(t+\frac{1}{2})^{1.1} + xy\alpha + \sqrt{t+\frac{1}{2}}z^3 + z(y-x)\alpha
    \end{array}
\right]$  \\
    \bottomrule
    \end{tabular}
    \label{table:manufact_sols:dimred_elastics}
\end{table}

\subsection{Modeling Nonlinear Elasticity}
\label{sec:elnl}
In Experiment~4, we consider a system with non-constant Young's modulus $E=E(\boldsymbol{\upvarepsilon})$. The stiffness of a material usually decreases with applied strain, so we choose a decreasing function
\begin{equation}
    E(\boldsymbol{\upvarepsilon}) = \frac{5}{\sqrt{20 + ||\boldsymbol{\upvarepsilon}||_{\textsc{f}}}},
\end{equation}
where $||\cdot||_{\textsc{f}}$ is the Frobenius norm. With this choice, the assumption of $E=1$ in the PBM is obviously not true, but still ``in the ball park'' for the manufactured solutions we consider here. These manufactured solutions, which have the same displacement $\boldsymbol{u}$ as the solutions discussed in Section~\ref{sec:el}, are presented in Table~\ref{table:manufact_sols:nl_elastic}.\footnote{Of course, since we use the same $\boldsymbol{u}$ but a different $E$, $\boldsymbol{\sigma}$ and $\boldsymbol{f}$ are different here than for Experiments~1 and~2. For the sake of brevity, we have not written out $\boldsymbol{\sigma}$ or $\boldsymbol{f}$ for any of out manufactured solutions here. However, they can be found in \citep{sorbo2022f}.}

The non-linearity greatly increases the cost of generating reference data. Therefore. only 500 time steps, $10\times10$ elements and 5 initializations of DDM and CoSTA were used in this experiment. The DNN learning rate was also increased to $8\cdot10^{-5}$.

\begin{table}
    \centering
    \caption{Manufactured solutions for nonlinear elastic problems.}
    \begin{tabular}{ ccc } 
    \toprule
    Label & $\boldsymbol{u}(t,x,y,\alpha)$ & $E$ \\ 
    \hline
    \vspace{8pt}
    \textit{n1} & $\left[
    \begin{array}{c}
        \sin(\pi(x+\alpha y))\cos(\alpha t) \\
        \cos(\pi(x+\alpha y))\sin(\alpha t)
    \end{array}
\right]$& $\frac{5}{\sqrt{20 + ||\boldsymbol{\upvarepsilon}||_{\textsc{f}}}}$\\ 
    \vspace{8pt}
    \textit{n2} & $\left[
    \begin{array}{c}
        \exp(\frac{(-tx^2 + y^2)}{(1+\alpha+t^2)}) \\
        \exp(\frac{(+tx^2 - y^2)}{(1+\alpha+t^2)})
    \end{array}
\right]$& $\frac{5}{\sqrt{20 + ||\boldsymbol{\upvarepsilon}||_{\textsc{f}}}}$\\ 
    \vspace{4pt}
    \textit{n3} & 
    $\left[
    \begin{array}{c}
        x^3 + y^2(t+0.5)^{1.5} + xy\alpha \\
        x^2 + y^3(t+0.5)^{1.1}+ xy\alpha
    \end{array}
\right]$& $\frac{5}{\sqrt{20 + ||\boldsymbol{\upvarepsilon}||_{\textsc{f}}}}$\\ 
    \bottomrule
    \end{tabular}
    \label{table:manufact_sols:nl_elastic}
\end{table}

\section{Results and discussion}
\label{sec:results}

For the exact implementation used for the experiments and figures in this work, see the first author's GitHub repository \citep{sorbo2022f}.

\subsection{Experiments 1 \& 2 -- Known and Unknown Load Term}
\label{sec:res:el}

\begin{figure*}
	\begin{subfigure}[b]{0.5\linewidth}
		\centering 
        \includegraphics[width=\textwidth]{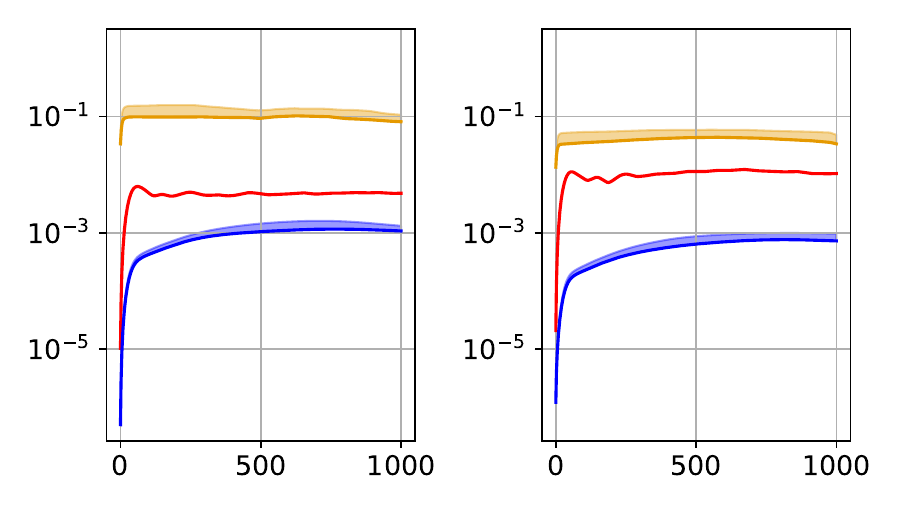}
        \subfigurecaptionInterpol{e1}
		\label{subfig:el_e_dev_i0}
	\end{subfigure}
	\begin{subfigure}[b]{0.5\linewidth}
		\centering 
        \includegraphics[width=\textwidth]{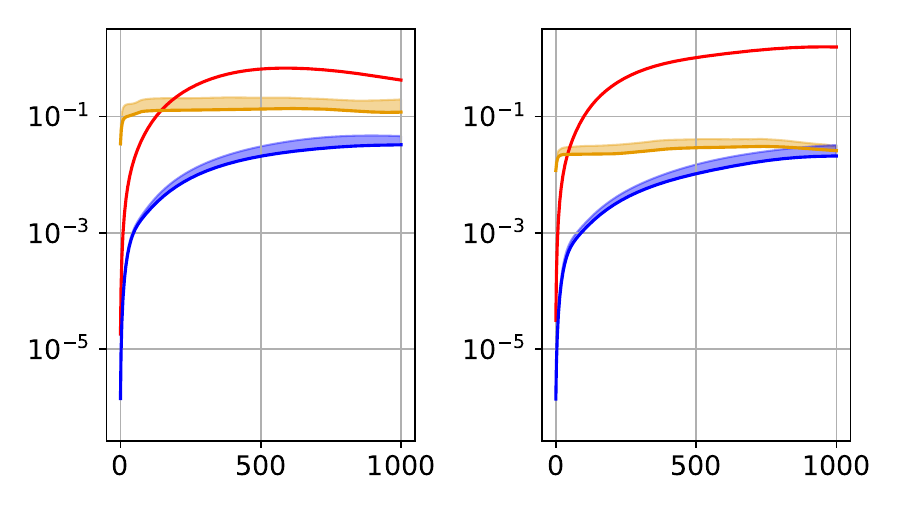}
        \subfigurecaptionInterpol{e1}
		\label{subfig:el_z_dev_i0}
	\end{subfigure}
	\\
	\begin{subfigure}[b]{0.5\linewidth}
		\centering 
        \includegraphics[width=\textwidth]{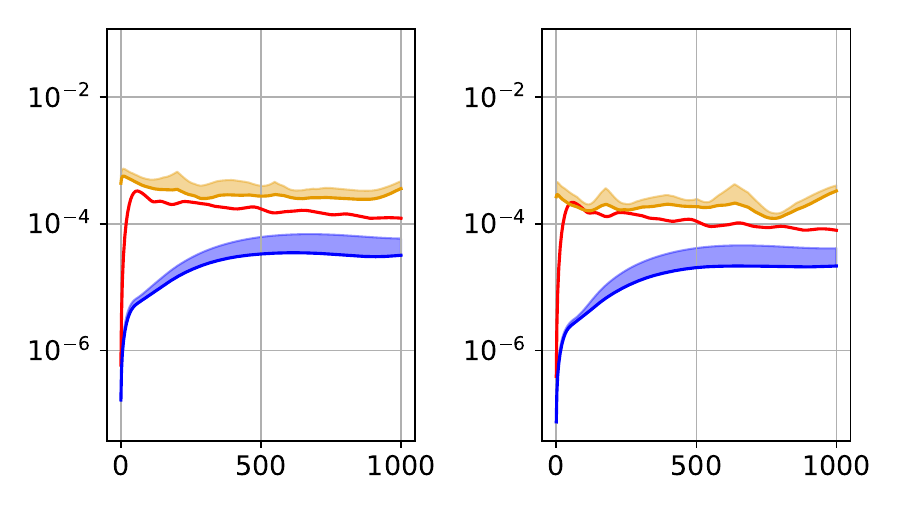}
        \subfigurecaptionInterpol{e2}
		\label{subfig:el_e_dev_i1}
	\end{subfigure}
	\begin{subfigure}[b]{0.5\linewidth}
		\centering 
        \includegraphics[width=\textwidth]{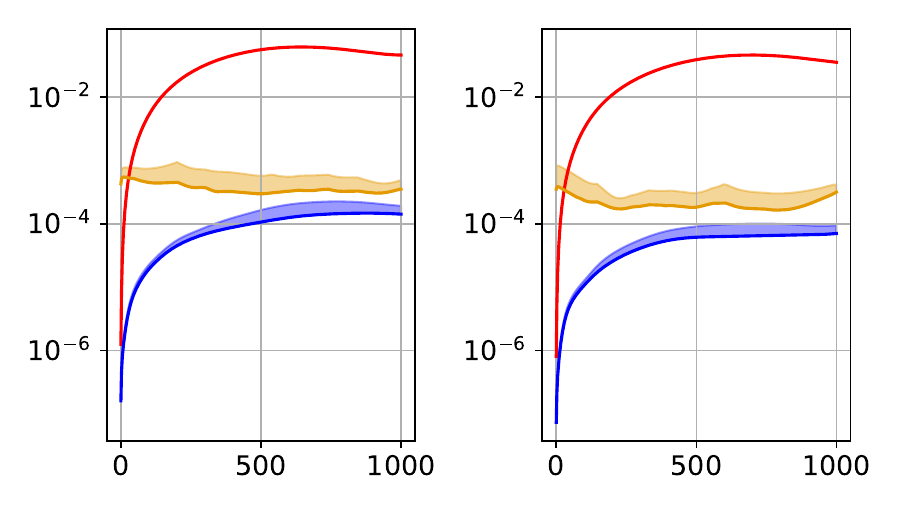}
        \subfigurecaptionInterpol{e2}
		\label{subfig:el_z_dev_i1}
	\end{subfigure}
	\\
	\begin{subfigure}[b]{0.5\linewidth}
		\centering 
        \includegraphics[width=\textwidth]{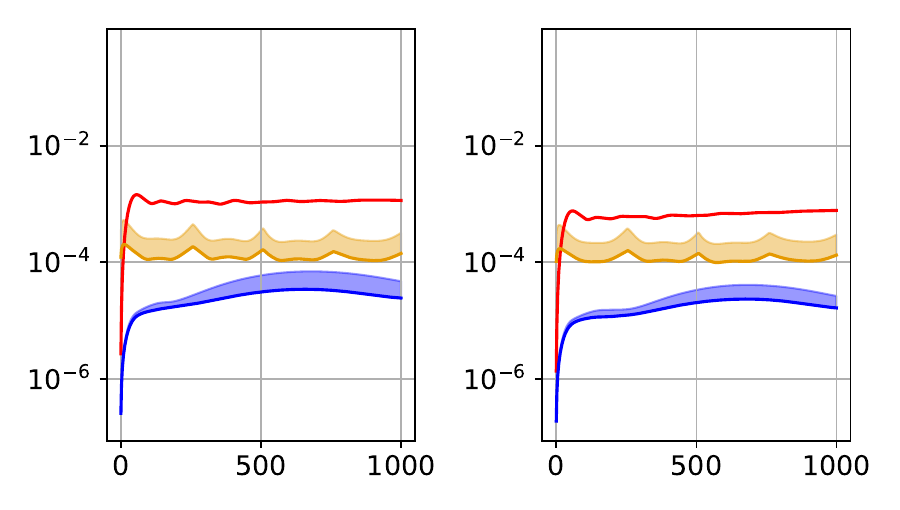}
        \subfigurecaptionInterpol{e3}
		\label{subfig:el_e_dev_i2}
	\end{subfigure}
	\begin{subfigure}[b]{0.5\linewidth}
		\centering 
        \includegraphics[width=\textwidth]{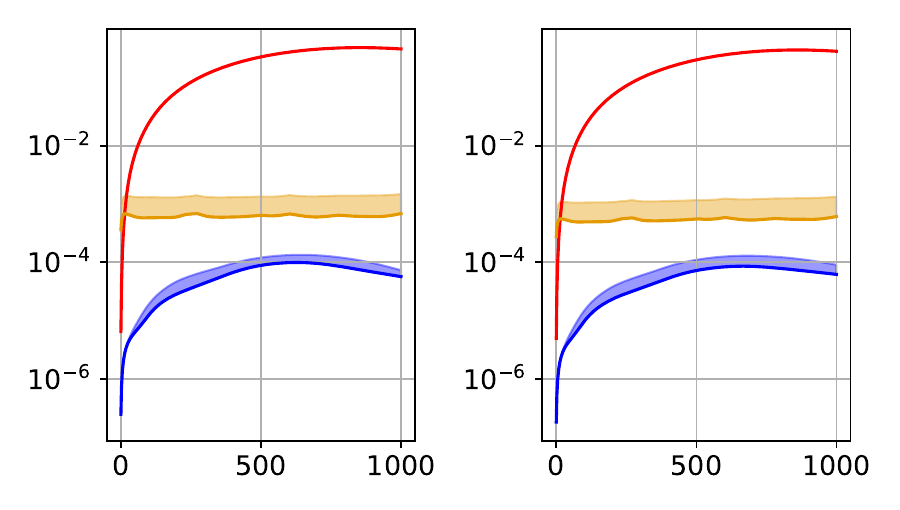}
        \subfigurecaptionInterpol{e3}
		\label{subfig:el_z_dev_i2}
		\centering 
	\end{subfigure}
\caption{Temporal development of relative $l_2$ error for solutions with correct source term (left) and zero source term (right) in interpolation scenarios. CoSTA is the most accurate method in all the cases.
($\pbmline$ PBM, $\ddmerr$ DDM, $\CoSTAerr$ CoSTA)
}
\label{fig:el_dev_i}
\end{figure*}

\begin{figure*}
	\begin{subfigure}[b]{0.5\linewidth}
		\centering 
        \includegraphics[width=\textwidth]{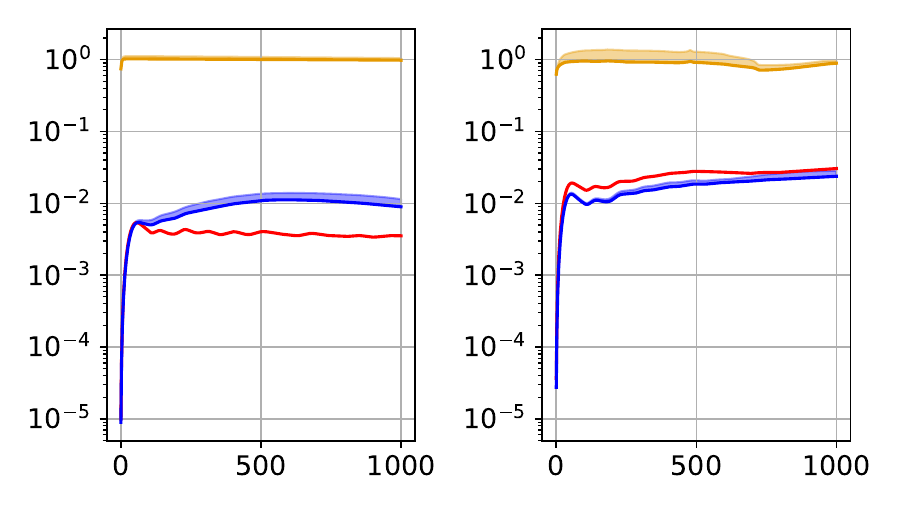}
        \subfigurecaptionExtrapol{e1}
		\label{subfig:el_f_dev_x0}
	\end{subfigure}
	\begin{subfigure}[b]{0.5\linewidth}
		\centering 
        \includegraphics[width=\textwidth]{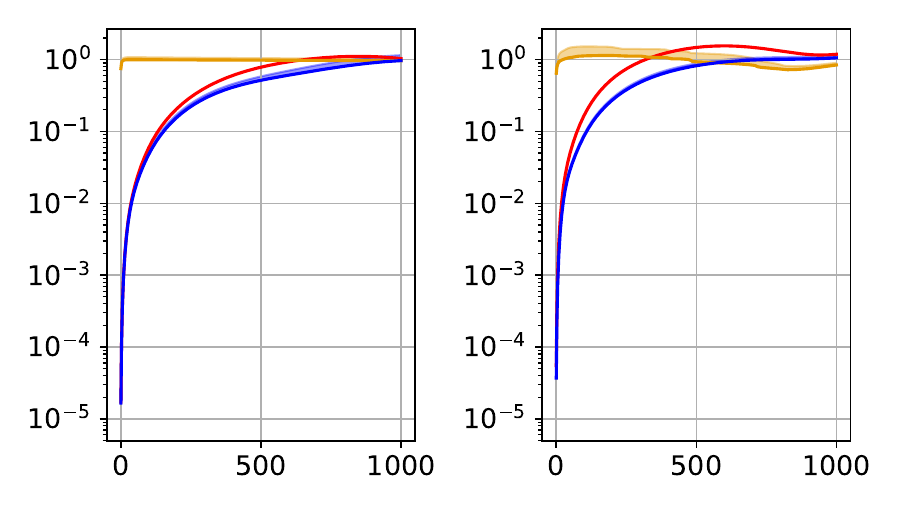}
        \subfigurecaptionExtrapol{e1}
		\label{subfig:el_z_dev_x0}
	\end{subfigure}
	\\
	\begin{subfigure}[b]{0.5\linewidth}
		\centering 
        \includegraphics[width=\textwidth]{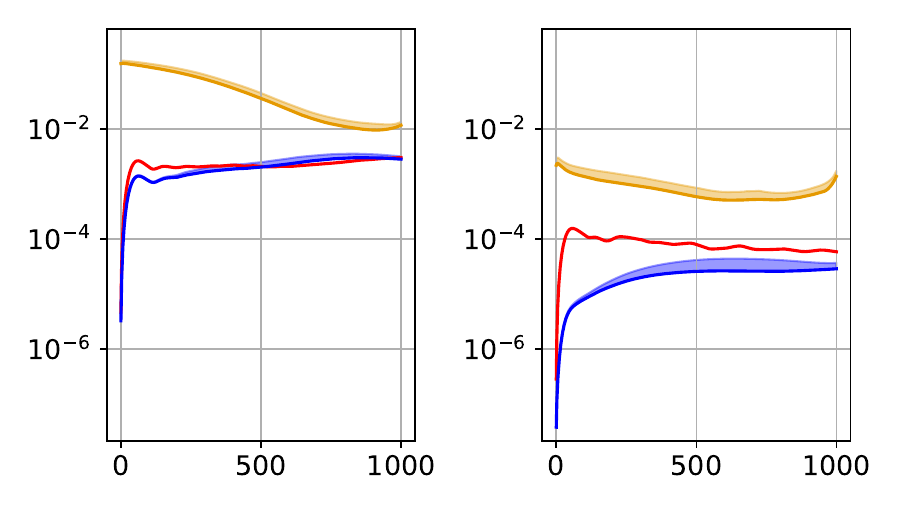}
        \subfigurecaptionExtrapol{e2}
		\label{subfig:el_f_dev_x1}
	\end{subfigure}
	\begin{subfigure}[b]{0.5\linewidth}
		\centering 
        \includegraphics[width=\textwidth]{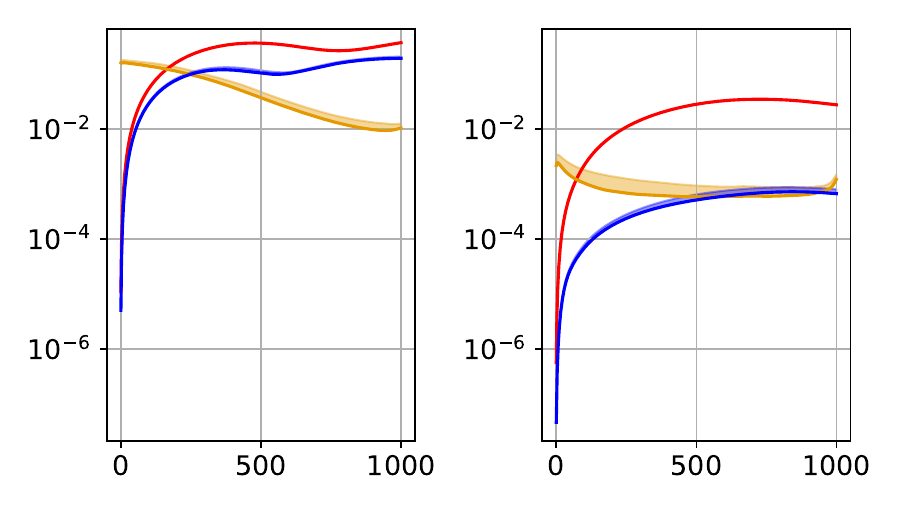}
        \subfigurecaptionExtrapol{e2}
		\label{subfig:el_z_dev_x1}
	\end{subfigure}
	\\
	\begin{subfigure}[b]{0.5\linewidth}
		\centering 
        \includegraphics[width=\textwidth]{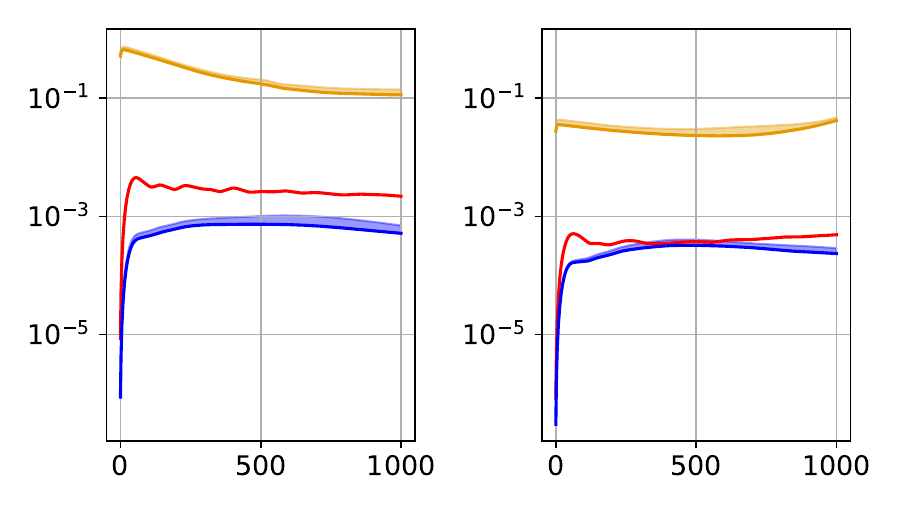}
        \subfigurecaptionExtrapol{e3}
		\label{subfig:el_f_dev_x2}
	\end{subfigure}
	\begin{subfigure}[b]{0.5\linewidth}
		\centering 
        \includegraphics[width=\textwidth]{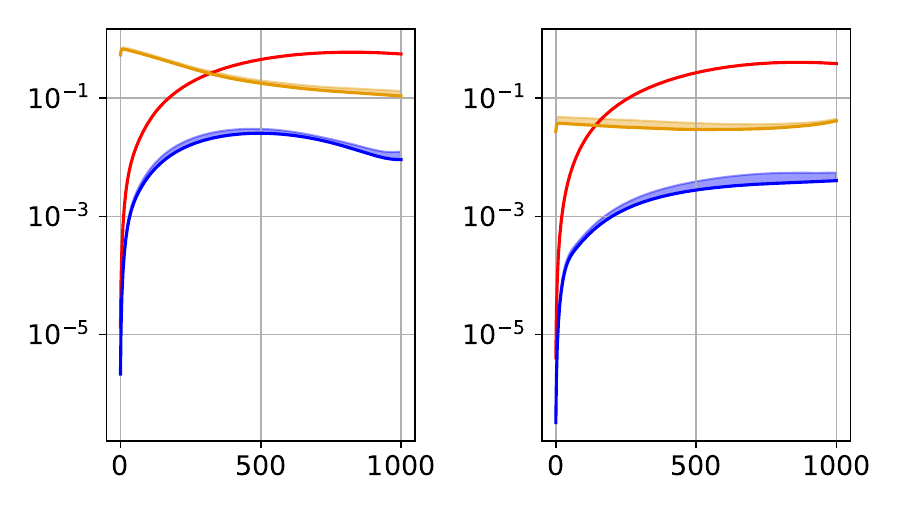}
        \subfigurecaptionExtrapol{e3}
		\label{subfig:el_z_dev_x2}
	\end{subfigure}
\caption{Temporal development of relative $l_2$ error for solutions with correct source term (left) and zero source term (right) in extrapolation scenarios. CoSTA is the best method, or among the best methods, in most cases, only being beaten once by each of the other methods.
($\pbmline$ PBM, $\ddmerr$ DDM, $\CoSTAerr$ CoSTA)
}
\label{fig:el_dev_x}
\end{figure*}

The temporal development of the models RRMSE is presented in Figure~\ref{fig:el_dev_i} for interpolation scenarios and Figure~\ref{fig:el_dev_x} for extrapolation scenarios. Results for the same manufactured solution with known and unknown load term are shown next to each other for easy comparison. We observe that CoSTA is more accurate than both PBM and DDM in all interpolation cases, and also in a significant majority of extrapolation cases. The only two cases where CoSTA is clearly not the most accurate model is $\alpha = -0.5$ for Solution~\textit{e1} with known source term and for Solution~\textit{e2} with unknown source term. In accordance with the discussion in Section~\ref{sec:introduction}, PBM performs at its best in the cases with no unknown physics, and DDM performs at its best in the interpolation scenarios.

The observation that CoSTA is generally more accurate than the PBM in Experiment~1 (with known load term) shows that CoSTA can be used to correct discretization error. This is in line with the findings in \citep{Blakseth2022dnn}, and suggests that CoSTA can be used to speed up expensive PBMs by permitting coarse discretizations without loss of accuracy. As for Experiment~2 (with unknown load term), since CoSTA outperforms DDM in most cases, it is evident that even a PBM severely affected by modeling error can be valuable in combination with a DNN.

From the top row of Figure~\ref{fig:el_dev_x}, it is clear that the purely data-driven model struggles in both extrapolation scenarios for Solution~\textit{e1}. Indeed, the DDM predictions have relative errors of roughly 100\%. Consequently, it is perhaps not so surprising that the data-driven component of CoSTA also struggles. Although the observed increase in error caused by CoSTA's DDM-generated correction term is undesirable, it is promising that CoSTA did not fail completely in a way similar to the DDM. This suggests that CoSTA is more robust than DDM.

For Solution~\textit{e2}, we observe that DDM performs quite well in all cases. This can perhaps be attributed to the simplicity of the solution, both in terms of spatial variability and its dependence on $\alpha$.\footnote{Solution~\textit{e2} at $t=0$ is largely dependent on $\alpha$, while its value at $t=1$ is much less so. Apart from the average slope, the general shape of the solution is also largely independent of $\alpha$.} Our results suggest that, in such simple cases, if the PBM is insufficiently accurate, it does not ease the learning task of the DDM, so pure DDM is more accurate. This is what we observe for Solution~\textit{e2} with $\alpha=-0.5$ and unknown load term.

Finally, we note that the error of CoSTA often fluctuates less and has a smaller standard deviation\footnote{Beware that, due to the logarithmic scaling of the $y$-axes, a larger shaded area does not necessarily imply a larger standard deviation (cf.\ e.g.\ the middle row, left half of Figure~\ref{fig:el_dev_i}).} than the DDM errors. This increased consistency contributes positively to the trustworthiness of CoSTA, as compared to DDM.

\subsection{Experiment 3 -- Dimensionality Reduction}
\label{sec:res:eldr}
\begin{figure*}
	\begin{subfigure}[b]{0.5\linewidth}
		\centering 
        \includegraphics[width=\textwidth]{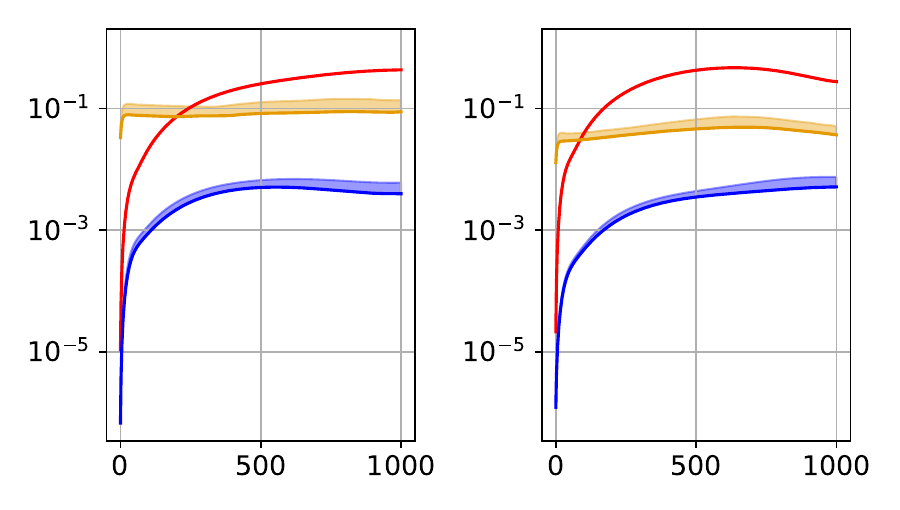}
        \subfigurecaptionInterpol{ed1}
		\label{subfig:el_dr_dev_i0}
	\end{subfigure}
	\begin{subfigure}[b]{0.5\linewidth}
		\centering 
        \includegraphics[width=\textwidth]{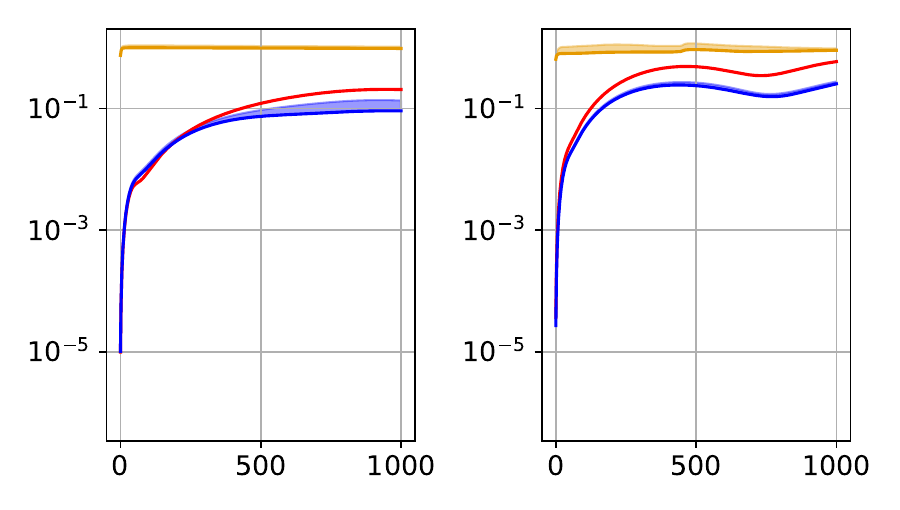}
        \subfigurecaptionExtrapol{ed1}
		\label{subfig:el_dr_dev_x0}
	\end{subfigure}
	\\
	\begin{subfigure}[b]{0.5\linewidth}
		\centering 
        \includegraphics[width=\textwidth]{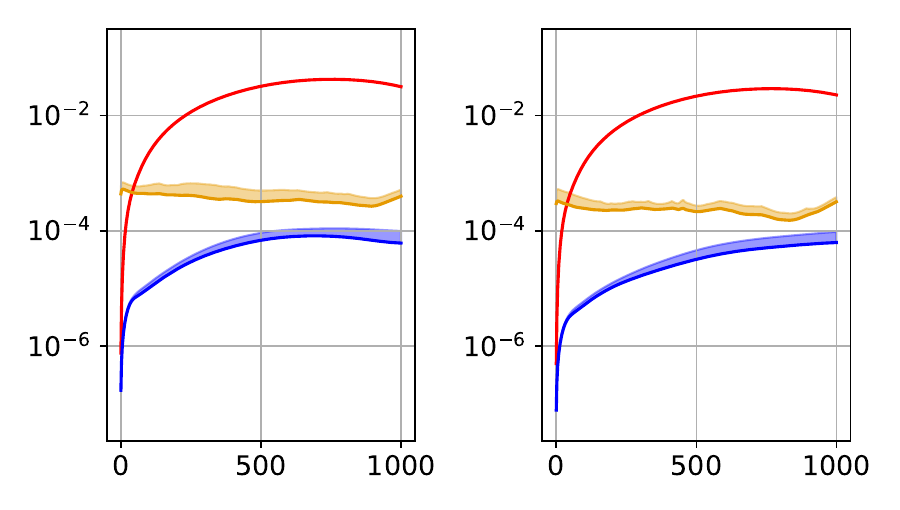}
        \subfigurecaptionInterpol{ed2}
		\label{subfig:el_dr_dev_i1}
	\end{subfigure}
	\begin{subfigure}[b]{0.5\linewidth}
		\centering 
        \includegraphics[width=\textwidth]{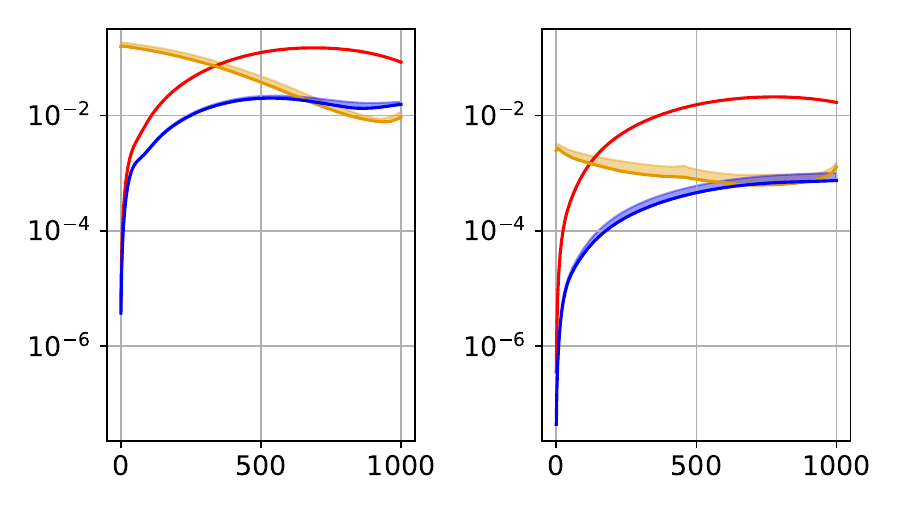}
        \subfigurecaptionExtrapol{ed2}
		\label{subfig:el_dr_dev_x1}
	\end{subfigure}
	\\
	\begin{subfigure}[b]{0.5\linewidth}
		\centering 
        \includegraphics[width=\textwidth]{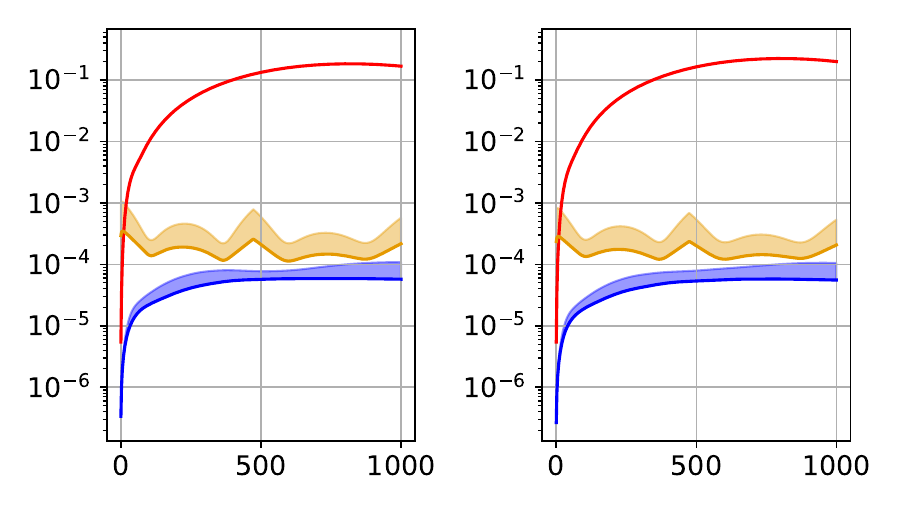}
        \subfigurecaptionInterpol{ed3}
		\label{subfig:el_dr_dev_i2}
	\end{subfigure}
	\begin{subfigure}[b]{0.5\linewidth}
		\centering 
        \includegraphics[width=\textwidth]{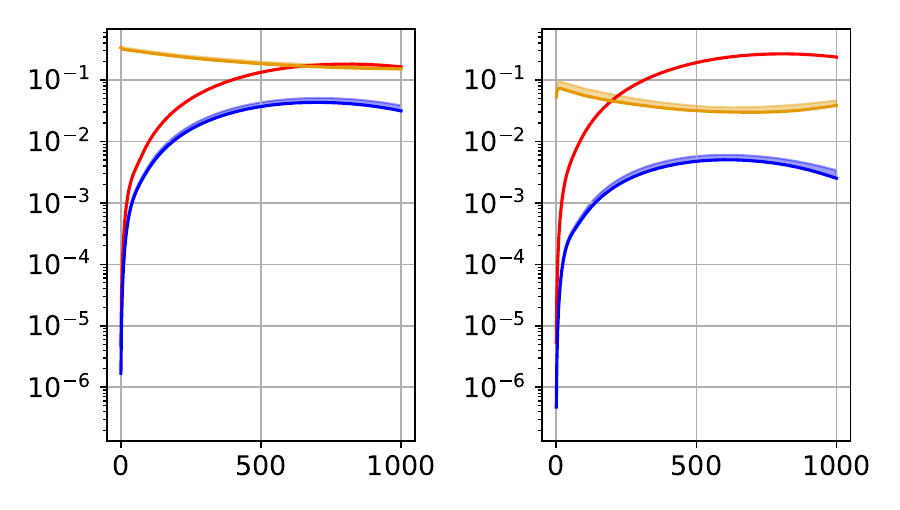}
        \subfigurecaptionExtrapol{ed3}
		\label{subfig:el_dr_dev_x2}
	\end{subfigure}
\caption{Temporal development of relative $l_2$ error for dimensional reduced linear elasticity in interpolation scenarios (left) and extrapolation scenarios (right).
Observe that CoSTA is more accurate than PBM in all cases, and better than DDM in all except $\alpha = -0.5$ for solution~\textit{ed2}.
($\pbmline$ PBM, $\ddmerr$ DDM, $\CoSTAerr$ CoSTA)
}
\label{fig:el_dr_dev}
\end{figure*}

The models' RRMSE for Experiment~3 are presented in Figure~\ref{fig:el_dr_dev}. For interpolation cases, CoSTA is consistently more accurate than the other two methods. DDM is in turn much better than PBM in these cases. For extrapolation, the DDM and PBM are alternately more accurate than each other, while CoSTA is clearly the most accurate model for Solutions~\textit{ed1} and~\textit{ed3}. For Solution~\textit{ed2} DDM is most accurate in the last time steps for $\alpha=-0.5$, while DDM and CoSTA are roughly equally accurate for $\alpha=2.5$. Overall, the results indicate that CoSTA is well-suited for correcting modeling error stemming from dimensionality reduction in PBMs.

\subsection{Experiment 4 -- Linearization of Nonlinear Elasticity}
\label{sec:res:elnl}

\begin{figure*}
	\begin{subfigure}[b]{0.5\linewidth}
		\centering 
        \includegraphics[width=\textwidth]{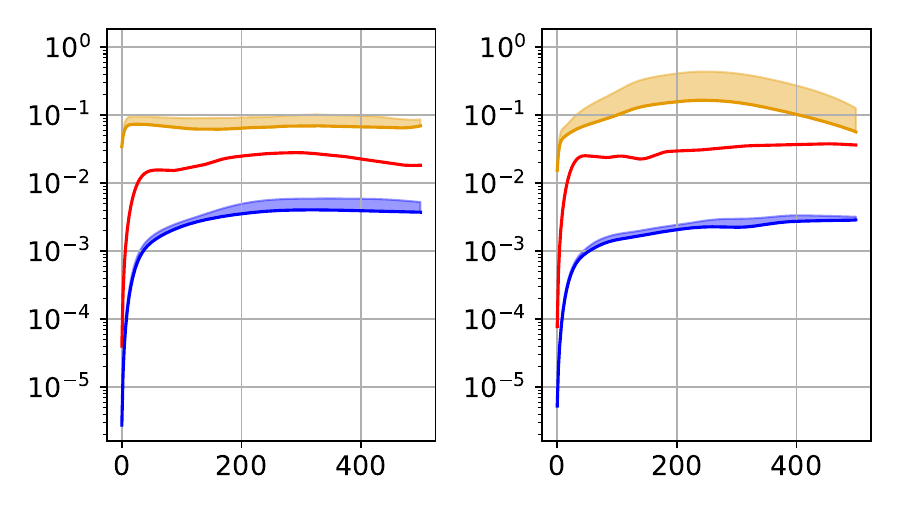}
        \subfigurecaptionInterpol{n1}
		\label{subfig:el_nl_dev_i0}
	\end{subfigure}
	\begin{subfigure}[b]{0.5\linewidth}
		\centering 
        \includegraphics[width=\textwidth]{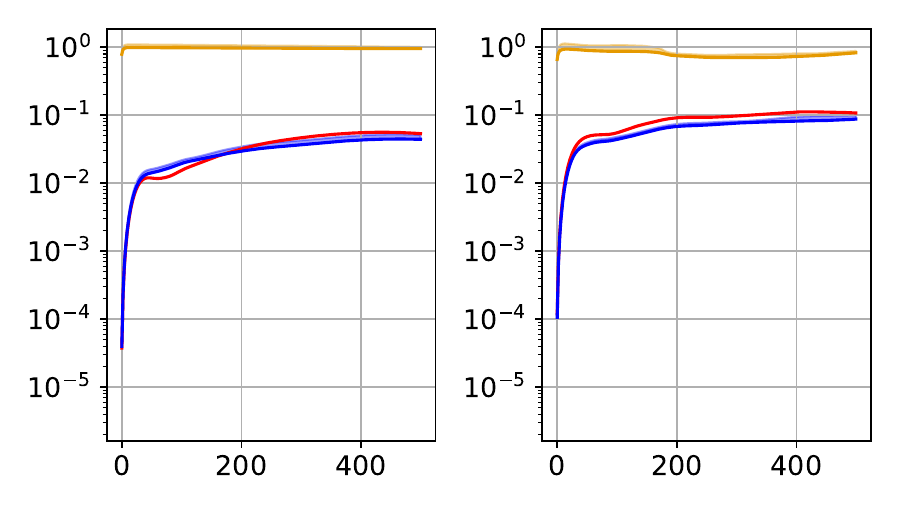}
        \subfigurecaptionExtrapol{n1}
		\label{subfig:el_nl_dev_x0}
	\end{subfigure}
	\\
	\begin{subfigure}[b]{0.5\linewidth}
		\centering 
        \includegraphics[width=\textwidth]{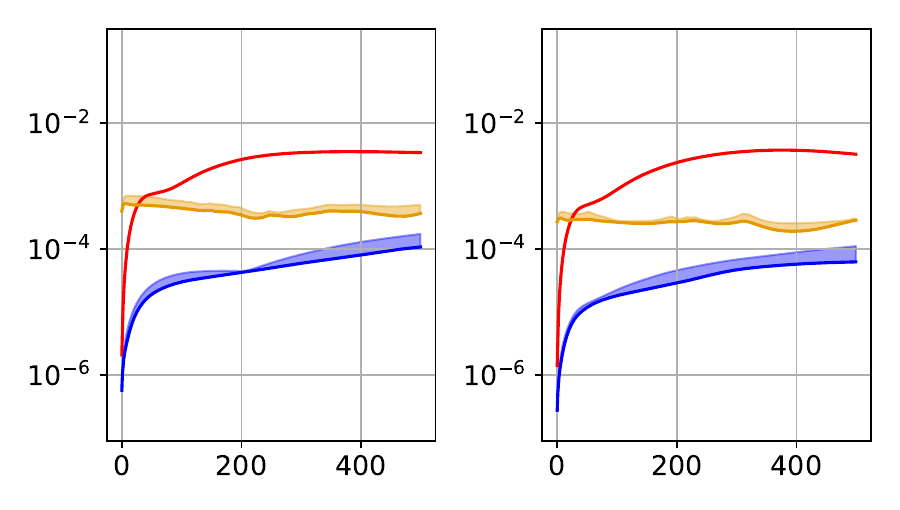}
        \subfigurecaptionInterpol{n2}
		\label{subfig:el_nl_dev_i1}
	\end{subfigure}
	\begin{subfigure}[b]{0.5\linewidth}
		\centering 
        \includegraphics[width=\textwidth]{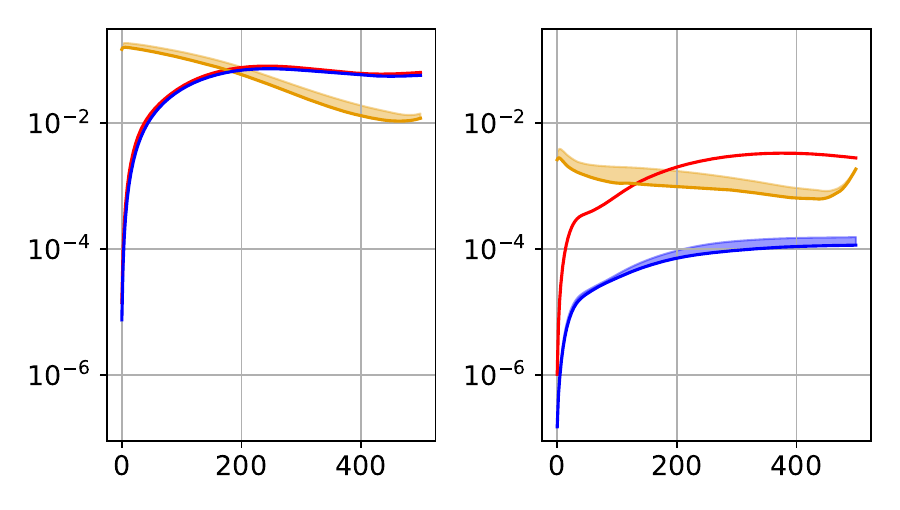}
        \subfigurecaptionExtrapol{n2}
		\label{subfig:el_nl_dev_x1}
	\end{subfigure}
	\\
	\begin{subfigure}[b]{0.5\linewidth}
		\centering 
        \includegraphics[width=\textwidth]{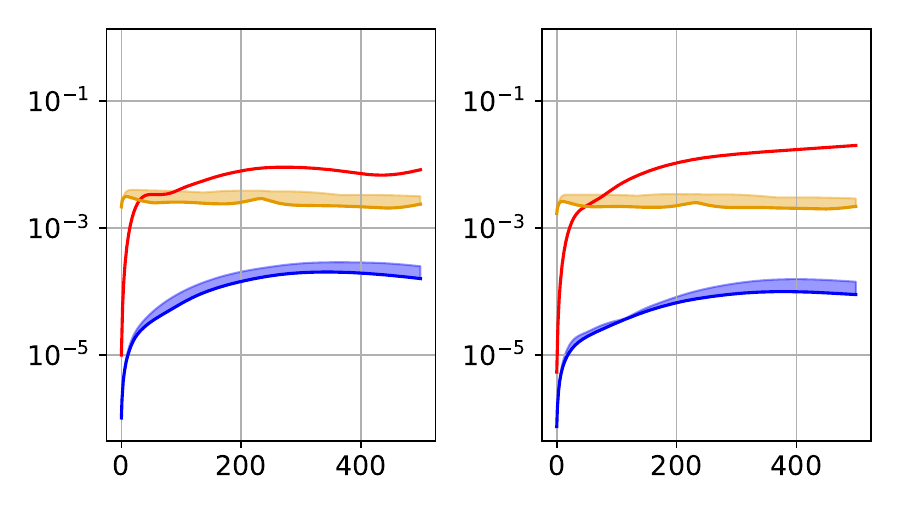}
        \subfigurecaptionInterpol{n3}
		\label{subfig:el_nl_dev_i2}
	\end{subfigure}
	\begin{subfigure}[b]{0.5\linewidth}
		\centering 
        \includegraphics[width=\textwidth]{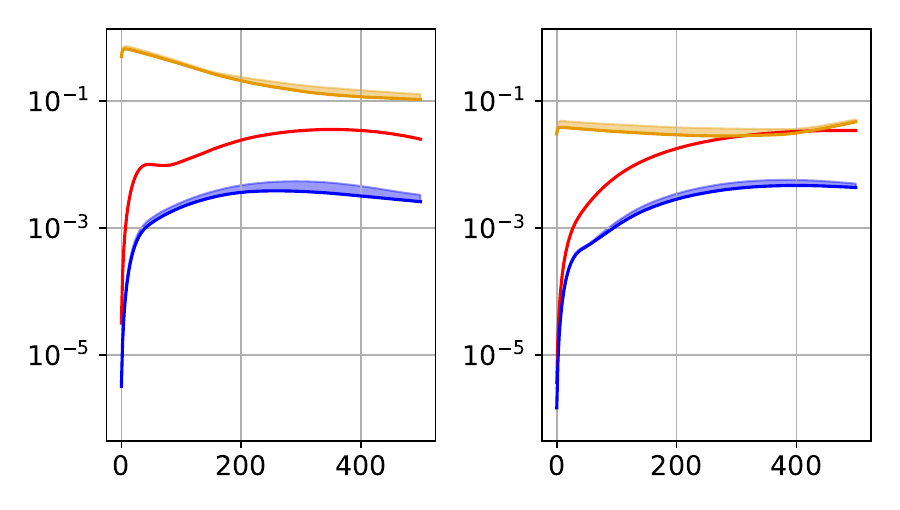}
        \subfigurecaptionExtrapol{n3}
		\label{subfig:el_nl_dev_x2}
	\end{subfigure}
\caption{Temporal development of relative $l_2$ error for nonlinear elasticity in interpolation scenarios (left) and extrapolation scenarios (right). While CoSTA is most accurate in all interpolation cases, the extrapolation results vary more.
($\pbmline$ PBM, $\ddmerr$ DDM, $\CoSTAerr$ CoSTA)
}
\label{fig:el_nl_dev}
\end{figure*}

The models' RRMSE for Experiment~4 are presented in Figure~\ref{fig:el_nl_dev}. CoSTA is clearly more accurate than PBM and DDM in all the interpolation cases. DDM is more accurate than PBM for all interpolation cases except Solution~\textit{n1}, where PBM is more accurate. 
In the extrapolation scenarios, CoSTA is significantly more accurate than PBM for Solution~\textit{n2} with $\alpha=2.5$ and Solution~\textit{n3}, while in the other cases the two model's are roughly equally accurate. 
DDM fails completely (~100\% RRMSE) for Solution~\textit{n1}, and is also the least accurate model for Solution~\textit{n3}. However, for Solution~\textit{n2}, DDM is more accurate than PBM in both cases and also more accurate than CoSTA for $\alpha=-0.5$. The latter case is the only one in this experiment where CoSTA is not the most accurate model.

\subsection{Result Comparison and Further Discussion}
\label{sec:more_discussion}

In this section, we summarize and provide further comments on the results of our numerical experiments.
We have discussed four experiments (concerning discretization error, unknown load term, dimensionality reduction and linearization) with 3 parametrized manufactured solutions and 4 $\alpha$-values combining to a total of 48 test scenarios.
In most scenarios, and in all interpolation scenarios, CoSTA has been the most accurate model. Still, there has been some variability in the relative performance of the three models.

\begin{figure}[tb]
    \centering
    \includegraphics[width=0.5\textwidth]{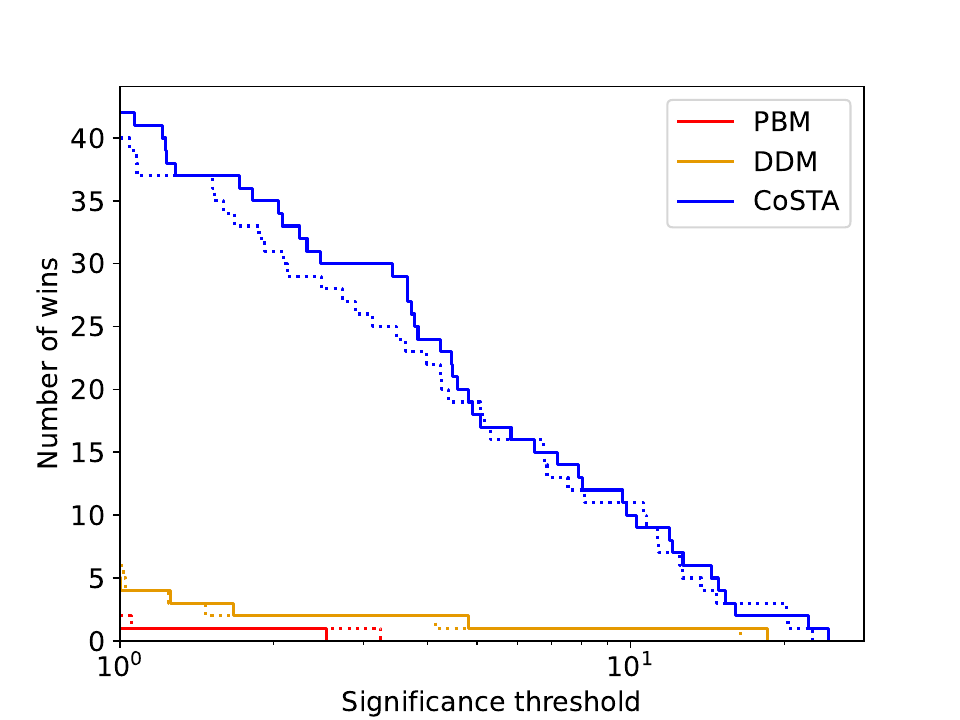}
    \caption{
    Overview of number of times each method was the most accurate at the final time step, at various degrees of significance. 
    For a value $x_1$ on the $x$-axis, the height of a curve is the number of times the method was better (had lower error) than both the other methods by a factor of at least $x_1$. 
    E.g. by observing the blue line at $10^1$, we see that in about 10 of the 48 tests, CoSTA won by a factor of at least $10^1$ (meaning both PBM and DDM had an error of more than $10$ times that of CoSTA).
    The solid line is the results of evaluating the means (of the CoSTA and DDM initializations), while the dotted line is from evaluating the mean + one standard deviation, as a way of penalizing inconsistency.}
    \label{fig:wins}
\end{figure}

\begin{figure}[tb]
    \centering
    \includegraphics[width=0.5\textwidth]{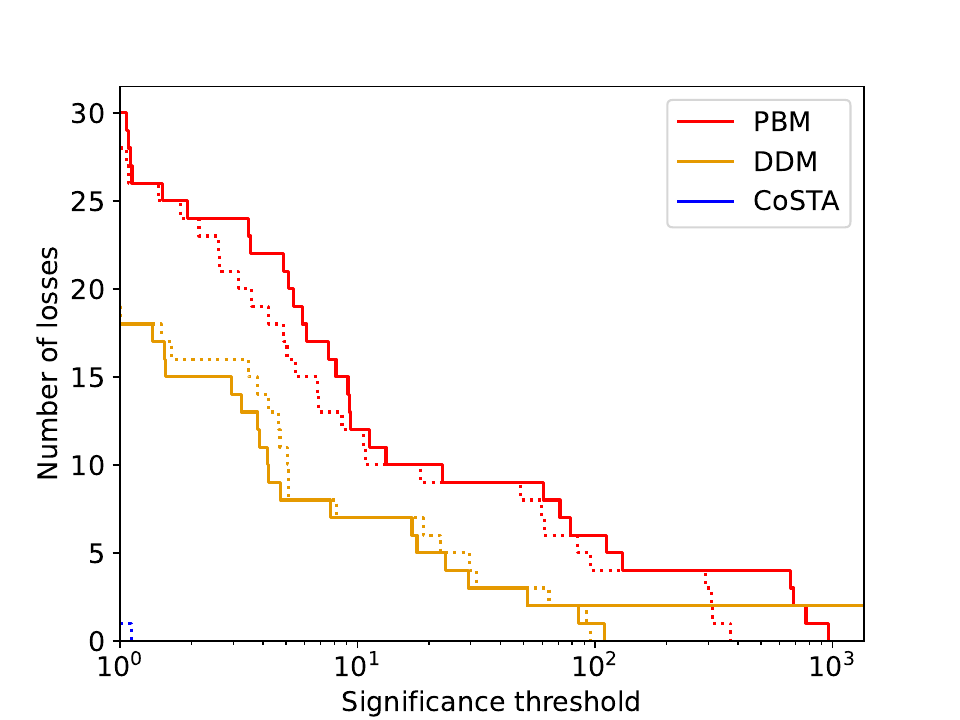}
    \caption{
    Overview of number of times each method was the \textit{least} accurate at the final time step, at various degrees of significance. 
    For a value $x_1$ on the $x$-axis, the height of a curve is the number of times a method was worse than both the other methods by a factor of at least $x_1$. 
    The solid line is the results of evaluating the means (of the CoSTA and DDM initializations), while the dotted line is from evaluating the mean + one standard deviation, as a way of penalizing inconsistency.
    We see that CoSTA was never the most inaccurate method, and barely the worst once when adding the standard deviation.}
    \label{fig:losses}
\end{figure}

Tables~\ref{tab:wins_by_alpha},~\ref{tab:wins_by_class} and~\ref{tab:wins_by_sol} show the number of times each model was the most accurate (i.e.\ had the lowest (mean) RRMSE at the final time step), categorized by $\alpha$-value, experiment and manufactured solution, respectively. Moreover, Figure~\ref{fig:wins} shows the number of times a model was more accurate than the other models by at least a factor $\delta$, as a function of significance threshold $\delta$. Similarly, Figure~\ref{fig:losses} shows the number of times each model was the least accurate, again as a function of the significance threshold. These tables and figures will form the basis for further discussion below.

First of all, Tables~\ref{tab:wins_by_alpha},~\ref{tab:wins_by_class} and~\ref{tab:wins_by_sol} and Figures~\ref{fig:wins} and~\ref{fig:losses} all support the statement that CoSTA is generally more accurate than PBM and DDM in our experiments.

In Table~\ref{tab:wins_by_alpha}, we observe that $\alpha=-0.5$ is the only $\alpha$-value for which CoSTA has not been the consistent winner. 
Given the discussion in Section~\ref{sec:introduction}, where PBM was labelled as more generalizable than DDM, it may seem strange that, in the extrapolation scenarios, the DDM is the most accurate model more frequently than the PBM. However, we believe that this can be explained by considering three factors, as explained in the following.

The first factor is that, out of the four scenarios with $\alpha=-0.5$ where the DDM is the most accurate model, three concerns the exponential manufactured solution. As touched upon briefly in Section~\ref{sec:res:el}, this solution has a qualitatively simple dependence on $\alpha$ which may be well suited for approximation by DNNs, although the underlying mechanism for this is unknown. The fourth case with $\alpha=-0.5$ where DDM ``wins'', is for Solution~\textit{e1} in Experiment~2. Here, although DDM is the most accurate, its RRMSE at the final time step is still close to 100\%, and its predictions consequently worthless from a practical standpoint. 

Secondly, across all 48 scenarios, the accuracy of DDM is more often than not higher than that of PBM (cf.\ Figures~\ref{fig:wins} and~\ref{fig:losses}). That is, on average, our DDM is more accurate than our PBM for the modeling problems considered in our experiments. Thus, DDM can afford to lose some accuracy due to poor generalizability and still be more accurate than PBM. Thirdly, the accuracy of DDM is worse in all extrapolation scenarios than in the corresponding interpolation scenarios. PBM does not exhibit a similar trend. Considering these factors, we believe that our results do not contradict the statement that PBM generalizes better than DDM, but rather support it.

\begin{table}[tb]
    \centering
    \caption{Overview of number times each method was the most accurate at final time step, for each value of $\alpha$ and in total.}
    \label{tab:wins_by_alpha}
    \begin{tabular}{ccccc}
        \toprule
        $\alpha$&PBM&DDM&CoSTA&Total  \\
        \hline
        -0.5&1&4&7&12\\
        0.7&0&0&12&12\\
        1.5&0&0&12&12\\
        2.5&0&1&11&12\\
        \hline
        Total&1&5&42&48\\
        \bottomrule
    \end{tabular}
\end{table}

\begin{table}[tb]
    \centering
    \caption{Overview of number times each method was the most accurate at final time step, for each type of extra modeling error.}
    \label{tab:wins_by_class}
    \begin{tabular}{ccccc}
        \toprule
        modeling error&PBM&DDM&CoSTA&Total  \\
        \hline
        None&1&0&11&12\\
        Zero source term&0&3&9&12\\
        Dimensionality reduction&0&1&11&12\\
        Linearization&0&1&11&12\\
        \bottomrule
    \end{tabular}
\end{table}

\begin{table}[tb]
    \centering
    \caption{Overview of number times each method was the most accurate at final time step, for each manufactured solution.}
    \label{tab:wins_by_sol}
    \begin{tabular}{ccccc}
        \toprule
        Solution&PBM&DDM&CoSTA&Total  \\
        \hline
        \textit{e1}/\textit{ed1}/\textit{n1} (sinusoidal) &1&2&13&16\\
        \textit{e2}/\textit{ed2}/\textit{n2} (exponential) &0&3&13&16\\
        \textit{e3}/\textit{ed3}/\textit{n3} (polynomial) &0&0&16&16\\
        \bottomrule
    \end{tabular}
\end{table}

Both DDM and CoSTA have, at times, exhibited noticeable variability in their predictions. (DDM more so than CoSTA.) Consistency is important for trustworthiness, and to avoid having to run several instances of the same model, and should therefore be taken into account when comparing the methods. One way of penalizing the variance when summarizing the results is to compare final mean error plus a standard deviation. The dotted lines in Figures \ref{fig:wins} and \ref{fig:losses} show the results of such a penalization. We see that CoSTA and DDM both perform a bit worse, and PBM a bit better, with this method of evaluation. But the difference between the dotted and solid lines is quite insignificant. In other words, the difference in accuracy between the methods, are generally much larger than the difference in accuracy between differently initialized instances of the same model.\footnote{This is easy to see in the error plots -- while the blue shaded area showing the standard deviation of CoSTA error at times is relatively thick, the distance to the other lines is usually considerably larger. 
Remember also that the thickness of the shaded area is also logarithmicaly scaled, so when the mean prediction is very accurate, a thick line does not indicate a big variance.
For example, in Figure~\ref{subfig:el_nl_dev_i2}, despite the variance, CoSTA is consistently much more accurate than the other methods.}

The above discussion on the scenarios where DDM is the most accurate model can shed light on the question of when to use CoSTA. It seems that, if the variable to predict is in some sense mathematically well-behaved in space and time, predicting a corrective source term is not necessarily simpler than predicting the solution directly. This appears to be especially true if 1) the qualitative behaviour of the system is similar across the relevant scenarios, 2) pure PBM offers inaccurate (or worse, qualitatively incorrect) predictions. In such scenarios, and for simpler solutions, DDM might outperform CoSTA.

On the opposite end of the spectrum, if a PBM of sufficient accuracy and acceptable computational cost is available, there is no need to use CoSTA. If the relevant scenarios are very different from those used to train CoSTA's DDM-component, CoSTA might even perform worse than pure PBM, as seen for Solution~\textit{e1} with $\alpha=-0.5$ in Experiment~1 (cf.\ Figure~\ref{fig:el_dev_x}).

We highlight that the possible limitations described above are rather weak. Indeed, in the vast majority of cases we have considered, the CoSTA model is significantly more accurate than the PBM and DDM models -- sometimes by more than one order of magnitude. Moreover, CoSTA's superior accuracy persists across a selection of different, commonly encountered error sources: discretization error, unknown physics, dimensionality reduction and linearization.
Overall, we conclude that CoSTA is the superior model in our numerical experiments.


\section{Conclusion and future work}
\label{sec:conclusionandfuturework}
In our research endeavors focused on elasticity modeling, we embarked on a series of numerical experiments aimed at showcasing the effectiveness of CoSTA in amalgamating physics-based modeling (PBM) and data-driven modeling (DDM) techniques into a comprehensive unified model. The results of these experiments, as detailed in this study, have revealed the considerable potential of CoSTA in bridging the gap between traditional PBM and DDM methodologies. More precisely:

\begin{itemize}
    \item When it comes to evaluating the accuracy of our CoSTA model, our findings present a compelling case. Across a total of 48 distinct experimental scenarios, our CoSTA model consistently outperformed its individual PBM and DDM counterparts in 42 of these cases. In particular, in numerous instances, the superiority of CoSTA was not marginal but rather substantial, surpassing its constituent models by more than one order of magnitude. This outstanding performance underscores the robustness and adaptability of CoSTA to address a wide range of engineering applications and error sources.
    \item The experiments carried out in this study encompassed a diverse set of error sources frequently encountered in engineering contexts. These included the discretization error, the error stemming from unknown physics, the dimensionality reduction error, and the linearization error. Our investigations reaffirmed earlier observations from related research \citep{Blakseth2022cpb, Blakseth2022dnn}, demonstrating that CoSTA excels at mitigating modeling errors attributed to discretization and errors related to the incorporation of unknown physics. Importantly, our work stands out as the first comprehensive demonstration of CoSTA's capabilities in correcting modeling errors arising from dimensionality reduction and linearization within PBMs.
    \item Crucially, our findings also highlight the robustness of CoSTA. Even when accounting for randomness introduced by deep neural network (DNN) initialization and the optimization process through stochastic gradient descent, CoSTA's superiority remained evident and persistent.
    \item Furthermore, it is worth noting a significant milestone achieved in this study. This work marks the first successful application of CoSTA in the context of PBMs based on the finite element method, showcasing the versatility of CoSTA across diverse modeling techniques. Furthermore, to foster transparency and promote further research in this area, we have made the implementation of our CoSTA model open source \citep{sorbo2022f}, thus allowing the broader scientific community to benefit from and build on our findings.
\end{itemize}
In summary, our research underscores the compelling potential of CoSTA as a unifying framework to combine PBM and DDM techniques, demonstrating its exceptional accuracy, adaptability, and robustness across a variety of engineering applications and error sources. However, we still foresee improvements to the CoSTA method in the following directions: 

\begin{itemize}
    \item Addressing the black-box nature through the use of symbolic regression: While CoSTA has shown promising results, addressing its black-box nature is essential for a deeper understanding and broader applicability. Symbolic regression techniques, such as genetic programming or neural-symbolic methods, can be explored to extract interpretable mathematical expressions from the DDM model. This would help uncover the underlying physics and relationships captured by CoSTA.
    \item Exploiting the correlation in time: To enhance CoSTA's ability to model transient systems and temporal dependencies, recurrent neural networks (RNNs) like Long Short-Term Memory (LSTM) or transformer-based architectures can be incorporated. These models excel at capturing sequential patterns and could be especially useful for systems with complex temporal dynamics.
    \item Dealing with noise: To make CoSTA more robust in practical scenarios, incorporating a denoising step into the pipeline is crucial. This step could involve training CoSTA on data with synthetic noise to improve its resilience to real-world noise. Additionally, exploring advanced denoising techniques, such as autoencoders can help enhance the model's performance in noisy environments.
    \item Using real-time measurements: To bridge the gap between laboratory experiments and real-world applications, it's essential to consider the use of real-time measurements. Building a latent space Reduced Order Model (ROM) that incorporates CoSTA as a component can facilitate real-time predictions and control. This would involve integrating CoSTA with online sensor data and dynamic system simulations to make it practical for real-world applications.
\end{itemize}

By addressing these four points, CoSTA's capabilities, interpretability, robustness, and real-world applicability can be further enhanced, ultimately advancing the field of computational modeling and simulation for various physical phenomena.

\section*{Acknowledgments}
This publication has been prepared as part of NorthWind (Norwegian Research Centre on Wind Energy) co-financed by the Research Council of Norway, industry and research partners. Read more at www.northwindresearch.no (grant no. 321954). O.S. gratefully acknowledges the National Science Foundation support (DMS-2012255).

\bibliographystyle{AR}
\bibliography{ref}

\begin{thebibliography}{10}
\expandafter\ifx\csname url\endcsname\relax
  \def\url#1{\texttt{#1}}\fi
\expandafter\ifx\csname urlprefix\endcsname\relax\def\urlprefix{URL }\fi
\expandafter\ifx\csname href\endcsname\relax
  \def\href#1#2{#2} \def\path#1{#1}\fi

\bibitem{Wei2018PredictingTE}
H.~Wei, S.~Zhao, Q.~Rong, H.~Bao,
  \href{https://www.sciencedirect.com/science/article/pii/S0017931018317423}{Predicting
  the effective thermal conductivities of composite materials and porous media
  by machine learning methods}, {\em International Journal of Heat and Mass
  Transfer}, 127:908--916 (2018).
\newblock

\bibitem{Baldi2016JetSC}
P.~Baldi, K.~T. Bauer, C.~Eng, P.~Sadowski, D.~Whiteson,
  \href{https://dx.doi.org/DOI: 10.1103/PhysRevD.93.094034}{Jet substructure
  classification in high-energy physics with deep neural networks}, {\em
  Physical Review D}, 93:094034 (2016).
\newblock

\bibitem{IbarraBerastegi2015ShorttermFO}
G.~Ibarra-Berastegi, J.~Saénz, G.~Esnaola, A.~Ezcurra, A.~Ulazia,
  \href{https://www.sciencedirect.com/science/article/pii/S0029801815002309}{Short-term
  forecasting of the wave energy flux: Analogues, random forests, and
  physics-based models}, {\em Ocean Engineering}, 104:530--539 (2015).
\newblock

\bibitem{Jia2021PhysicsGuidedML}
X.~Jia, J.~D. Willard, A.~Karpatne, J.~S. Read, J.~A. Zwart, M.~S. Steinbach,
  V.~Kumar,
  \href{https://dx.doi.org/https://doi.org/10.1145/3447814}{Physics-guided
  machine learning for scientific discovery: An application in simulating lake
  temperature profiles}, {\em ACM/IMS Transactions on Data Science}, 2:1 -- 26
  (2021).
\newblock

\bibitem{san2021hybrid}
O.~San, A.~Rasheed, T.~Kvamsdal,
  \href{https://dx.doi.org/https://doi.org/10.1002/gamm.202100007}{Hybrid
  analysis and modeling, eclecticism, and multifidelity computing toward
  digital twin revolution}, {\em GAMM-Mitteilungen}, 44:e202100007 (2021).
\newblock

\bibitem{8972429}
A.~Rasheed, O.~San, T.~Kvamsdal,
  \href{https://dx.doi.org/10.1109/ACCESS.2020.2970143}{Digital twin: Values,
  challenges and enablers from a modeling perspective}, {\em IEEE Access},
  8:21980--22012 (2020).
\newblock

\bibitem{willard2020integrating}
J.~Willard, X.~Jia, S.~Xu, M.~Steinbach, V.~Kumar, Integrating physics-based
  modeling with machine learning: A survey, {\em arXiv preprint
  arXiv:2003.04919} (2020).

\bibitem{amos_optnet_2017}
B.~Amos, J.~Z. Kolter,
  \href{https://proceedings.mlr.press/v70/amos17a.html}{{OptNet}:
  {Differentiable} {Optimization} as a {Layer} in {Neural} {Networks}}, in {\em
  International {Conference} on {Machine} {Learning}}, PMLR, 2017, pp.
  136--145.

\bibitem{avilabelbuteperes_end_2018}
F.~de~Avila Belbute-Peres, K.~Smith, K.~Allen, J.~Tenenbaum, J.~Z. Kolter,
  \href{https://proceedings.neurips.cc/paper/2018/file/842424a1d0595b76ec4fa03c46e8d755-Paper.pdf}{End-to-end
  differentiable physics for learning and control}, in S.~Bengio, H.~Wallach,
  H.~Larochelle, K.~Grauman, N.~Cesa-Bianchi, R.~Garnett (eds.), {\em Advances
  in Neural Information Processing Systems}, vol.~31, Curran Associates, Inc.,
  2018, p.~1.

\bibitem{yu2020sds}
Y.~Yu, H.~Yao, Y.~Liu,
  \href{https://www.sciencedirect.com/science/article/pii/S0952197620302670}{Structural
  dynamics simulation using a novel physics-guided machine learning method},
  {\em Engineering Applications of Artificial Intelligence}, 96:103947 (2020).
\newblock

\bibitem{quarteroni2014reduced}
A.~Quarteroni, G.~Rozza, {\em Reduced order methods for modeling and
  computational reduction}, vol.~9, Springer, New York, 2014.

\bibitem{fonn2019dcp}
E.~Fonn, H.~v. Brummelen, T.~Kvamsdal, A.~Rasheed,
  \href{https://dx.doi.org/https://doi.org/10.1016/j.cma.2018.11.038}{Fast
  divergence-conforming reduced basis methods for steady navier–stokes flow},
  {\em Computer Methods in Applied Mechanics and Engineering}, 346:486--512
  (2019).
\newblock

\bibitem{pawar2019aet}
S.~{Pawar}, S.~E. {Ahmed}, O.~{San}, A.~{Rasheed},
  \href{https://dx.doi.org/https://doi.org/10.3390/math8040570}{An
  evolve-then-correct reduced order model for hidden fluid dynamics}, {\em
  Mathematics}, 8:570 (2020).
\newblock

\bibitem{pawar2020ddr}
S.~Pawar, S.~E. Ahmed, O.~San, A.~Rasheed,
  \href{https://dx.doi.org/https://doi.org/10.1063/5.0002051}{Data-driven
  recovery of hidden physics in reduced order modeling of fluid flows}, {\em
  Physics of Fluids}, 32:036602 (2020).
\newblock

\bibitem{raissi2019physics}
M.~Raissi, P.~Perdikaris, G.~E. Karniadakis,
  \href{https://dx.doi.org/https://doi.org/10.1016/j.jcp.2018.10.045}{Physics-informed
  neural networks: A deep learning framework for solving forward and inverse
  problems involving nonlinear partial differential equations}, {\em Journal of
  Computational Physics}, 378:686--707 (2019).
\newblock

\bibitem{zobeiry_physics_2021}
N.~Zobeiry, K.~D. Humfeld,
  \href{https://www.sciencedirect.com/science/article/pii/S0952197621000798}{A
  physics-informed machine learning approach for solving heat transfer equation
  in advanced manufacturing and engineering applications}, {\em Engineering
  Applications of Artificial Intelligence}, 101:104232 (2021).
\newblock

\bibitem{arnold_state_2021}
F.~Arnold, R.~King,
  \href{https://www.sciencedirect.com/science/article/pii/S0952197621000427}{{State-space
  modeling for control based on physics-informed neural networks}}, {\em
  Engineering Applications of Artificial Intelligence}, 101:104195 (2021).
\newblock

\bibitem{shen_physics_2021}
S.~Shen, H.~Lu, M.~Sadoughi, C.~Hu, V.~Nemani, A.~Thelen, K.~Webster, M.~Darr,
  J.~Sidon, S.~Kenny,
  \href{https://www.sciencedirect.com/science/article/pii/S0952197621001421}{A
  physics-informed deep learning approach for bearing fault detection}, {\em
  Engineering Applications of Artificial Intelligence}, 103:104295 (2021).
\newblock

\bibitem{BILLAH2023106336}
M.~M. Billah, A.~I. Khan, J.~Liu, P.~Dutta,
  \href{https://www.sciencedirect.com/science/article/pii/S2352492823010279}{Physics-informed
  deep neural network for inverse heat transfer problems in materials}, {\em
  Materials Today Communications}, 35:106336 (2023).
\newblock

\bibitem{SUN2023103525}
Z.~Sun, H.~Du, C.~Miao, Q.~Hou,
  \href{https://www.sciencedirect.com/science/article/pii/S0965997823001163}{A
  physics-informed neural network based simulation tool for reacting flow with
  multicomponent reactants}, {\em Advances in Engineering Software}, 185:103525
  (2023).
\newblock

\bibitem{HALAMKA2023109351}
J.~Halamka, M.~Bartošák, M.~Španiel,
  \href{https://www.sciencedirect.com/science/article/pii/S0013794423003090}{Using
  hybrid physics-informed neural networks to predict lifetime under multiaxial
  fatigue loading}, {\em Engineering Fracture Mechanics}, 289:109351 (2023).
\newblock

\bibitem{Champion22445}
K.~Champion, B.~Lusch, J.~N. Kutz, S.~L. Brunton,
  \href{https://www.pnas.org/content/116/45/22445}{Data-driven discovery of
  coordinates and governing equations}, {\em Proceedings of the National
  Academy of Sciences}, 116:22445--22451 (2019).
\newblock

\bibitem{vaddireddy2020fes}
H.~Vaddireddy, A.~Rasheed, A.~E. Staples, O.~San,
  \href{https://dx.doi.org/https://doi.org/10.1063/1.5136351}{Feature
  engineering and symbolic regression methods for detecting hidden physics from
  sparse sensors}, {\em Physics of Fluids, Editor's pick}, 32:015113 (2020).
\newblock

\bibitem{ZHANG2023112349}
H.~Zhang, X.~Liu, G.~Zhang, Y.~Zhu, S.~Li, Q.~Qian, D.~Dai, R.~Che, T.~Xu,
  \href{https://www.sciencedirect.com/science/article/pii/S0927025623003439}{Deriving
  equation from data via knowledge discovery and machine learning: A study of
  young’s modulus of ti-nb alloys}, {\em Computational Materials Science},
  228:112349 (2023).
\newblock

\bibitem{MEYER2023105416}
K.~A. Meyer, F.~Ekre,
  \href{https://www.sciencedirect.com/science/article/pii/S002250962300220X}{Thermodynamically
  consistent neural network plasticity modeling and discovery of evolution
  laws}, {\em Journal of the Mechanics and Physics of Solids}, 180:105416
  (2023).
\newblock

\bibitem{pawar2021pgml}
S.~Pawar, O.~San, B.~Aksoylu, A.~Rasheed, T.~Kvamsdal,
  \href{https://dx.doi.org/https://doi.org/10.1063/5.0038929}{Physics guided
  machine learning using simplified theories}, {\em Physics of Fluids},
  33:011701 (2021).
\newblock

\bibitem{Robinson2022pgn}
H.~Robinson, S.~Pawar, A.~Rasheed, O.~San,
  \href{https://www.sciencedirect.com/science/article/pii/S0893608022002854}{Physics
  guided neural networks for modelling of non-linear dynamics}, {\em Neural
  Networks}, 154:333--345 (2022).
\newblock

\bibitem{pineda2022theseus}
L.~Pineda, T.~Fan, M.~Monge, S.~Venkataraman, P.~Sodhi, R.~T. Chen, J.~Ortiz,
  D.~DeTone, A.~Wang, S.~Anderson, J.~Dong, B.~Amos, M.~Mukadam, {Theseus: A
  Library for Differentiable Nonlinear Optimization}, {\em Advances in Neural
  Information Processing Systems} (2022).

\bibitem{Blakseth2022dnn}
S.~S. Blakseth, A.~Rasheed, T.~Kvamsdal, O.~San,
  \href{https://www.sciencedirect.com/science/article/pii/S0893608021004494}{Deep
  neural network enabled corrective source term approach to hybrid analysis and
  modeling}, {\em Neural Networks}, 146:181--199 (2022).
\newblock

\bibitem{maulik2019sgm}
R.~Maulik, O.~San, A.~Rasheed, P.~Vedula,
  \href{https://dx.doi.org/https://doi.org/10.1017/jfm.2018.770}{Sub-grid
  modelling for two-dimensional turbulence using neural networks}, {\em Journal
  of Fluid Mechanics}, 858:122--144 (2019).
\newblock

\bibitem{pawar2020apa}
S.~Pawar, O.~San, A.~Rasheed, P.~Vedula,
  \href{https://dx.doi.org/https://doi.org/10.1007/s00162-019-00512-z}{A priori
  analysis on deep learning of subgrid-scale parameterizations for {Kraichnan}
  turbulence}, {\em Theoretical and Computational Fluid Dynamics}, 34:429--455
  (2020).
\newblock

\bibitem{robinson2022novel}
H.~Robinson, E.~Lundby, A.~Rasheed, J.~T. Gravdahl,
  \href{https://www.sciencedirect.com/science/article/pii/S0952197623008072}{Deep
  learning assisted physics-based modeling of aluminum extraction process},
  {\em Engineering Applications of Artificial Intelligence}, 125:106623 (2023).
\newblock

\bibitem{blakseth2021ica}
S.~S. Blakseth, Introducing {CoSTA}: A deep neural network enabled approach to
  improving physics-based numerical simulations, Master's thesis, NTNU (2021).

\bibitem{Blakseth2022cpb}
S.~S. Blakseth, A.~Rasheed, T.~Kvamsdal, O.~San,
  \href{https://dx.doi.org/10.1016/j.asoc.2022.109533}{Combining physics-based
  and data-driven techniques for reliable hybrid analysis and modeling using
  the corrective source term approach}, {\em Applied Soft Computing},
  128:109533 (2022).
\newblock

\bibitem{linear_elasticity}
W.~S. Slaughter, {\em The Linearized Theory of Elasticity}, Birkhäuser,
  Boston, MA, 2002.

\bibitem{irgens2008continuum}
F.~Irgens, \href{https://books.google.no/books?id=q5dB7Gf4bIoC}{{\em Continuum
  Mechanics}}, Springer Berlin Heidelberg, 2008.

\bibitem{quarteroni2017nummod}
A.~Quarteroni, {\em Numerical models for differential problems}, 3rd edition,
  Springer International Publishing, 2017.

\bibitem{brenner2008FEM}
S.~C. Brenner, L.~R. Scott, {\em The mathematical theory of finite element
  methods}, 3rd edition, Springer, New York, 2008.

\bibitem{Courant1943VariationalMF}
R.~Courant, Variational methods for the solution of problems of equilibrium and
  vibrations, {\em Bulletin of the American Mathematical Society}, 49:1--23
  (1943).

\bibitem{Cybenko1989ApproximationBS}
G.~V. Cybenko,
  \href{https://dx.doi.org/https://doi.org/10.1007/BF02551274}{Approximation by
  superpositions of a sigmoidal function}, {\em Mathematics of Control, Signals
  and Systems}, 2:303--314 (1989).
\newblock

\bibitem{Hornik1989MultilayerFN}
K.~Hornik, M.~B. Stinchcombe, H.~L. White, Multilayer feedforward networks are
  universal approximators, {\em Neural Networks}, 2:359--366 (1989).

\bibitem{LESHNO1993861}
M.~Leshno, V.~Y. Lin, A.~Pinkus, S.~Schocken,
  \href{https://www.sciencedirect.com/science/article/pii/S0893608005801315}{Multilayer
  feedforward networks with a nonpolynomial activation function can approximate
  any function}, {\em Neural Networks}, 6:861--867 (1993).
\newblock

\bibitem{sym10110648}
I.~Nusrat, S.-B. Jang, \href{https://www.mdpi.com/2073-8994/10/11/648}{A
  comparison of regularization techniques in deep neural networks}, {\em
  Symmetry}, 10 (2018).
\newblock

\bibitem{Roache2001cvb}
P.~J. Roache, \href{https://dx.doi.org/https://doi.org/10.1115/1.1436090}{{Code
  Verification by the Method of Manufactured Solutions }}, {\em Journal of
  Fluids Engineering}, 124:4--10 (2001).
\newblock

\bibitem{tensorflow2015-whitepaper}
M.~Abadi, A.~Agarwal, P.~Barham, E.~Brevdo, Z.~Chen, C.~Citro, G.~S. Corrado,
  A.~Davis, J.~Dean, M.~Devin, S.~Ghemawat, I.~Goodfellow, A.~Harp, G.~Irving,
  M.~Isard, Y.~Jia, R.~Jozefowicz, L.~Kaiser, M.~Kudlur, J.~Levenberg,
  D.~Man\'{e}, R.~Monga, S.~Moore, D.~Murray, C.~Olah, M.~Schuster, J.~Shlens,
  B.~Steiner, I.~Sutskever, K.~Talwar, P.~Tucker, V.~Vanhoucke, V.~Vasudevan,
  F.~Vi\'{e}gas, O.~Vinyals, P.~Warden, M.~Wattenberg, M.~Wicke, Y.~Yu,
  X.~Zheng, \href{https://www.tensorflow.org/}{{TensorFlow}: Large-scale
  machine learning on heterogeneous systems}, Software available from
  tensorflow.org (2015).

\bibitem{lrelu}
A.~L. Maas, A.~Y. Hannum, A.~Y. Ng, Rectifier nonlinearities improve neural
  network acoustic models, in {\em Proceedings of the 30$^{th}$ International
  Conference on Machine Learning}, vol.~28, 2013, p.~3.

\bibitem{kingma2014aam}
D.~P. Kingma, J.~Ba, Adam: A method for stochastic optimization, {\em arXiv
  preprint arXiv:1412.6980} (2014).

\bibitem{sorbo2022f}
S.~Sørbø, {FEM\_CoSTA}, GitHub Repository
  \url{https://github.com/sondsorb/FEM_CoSTA} (2022).

\end{thebibliography}
\end{document}